\documentclass[a4paper,10pt]{amsart}

\usepackage[english]{babel}
\usepackage[english]{translator}
\usepackage[T1]{fontenc}
\usepackage[utf8]{inputenc}

\usepackage{amsmath}
\usepackage{amssymb}
\usepackage{amsfonts}
\usepackage{amstext}
\usepackage{amsthm}
\usepackage{bbm}

\usepackage{graphicx}
\usepackage{color}
\usepackage{psfrag}
\usepackage{subfigure}
\usepackage{float}

\usepackage{marvosym}

\usepackage{enumitem}
\usepackage{todonotes}
\usepackage{titlesec}

\usepackage{hyperref}

\usepackage[
  hmarginratio={1:1},     
  vmarginratio={1:1},     
  textwidth=460pt,        
  heightrounded,          
  bottom=3.8cm,
  top=3.8cm,
]{geometry}

\hypersetup{
  pdffitwindow=false,
  pdfhighlight=/O,
  pdfnewwindow,
  colorlinks=true,
  citecolor=red,            
  linkcolor=blue,             
  menucolor=blue,            
  urlcolor=blue,             
  pdfpagemode=UseOutlines,
  bookmarksnumbered=true,
  linktocpage,
  pdftitle={Freezing Traveling and Rotating Waves in Second Order Evolution Equations},
  pdfsubject={Freezing Traveling and Rotating Waves in Second Order Evolution Equations},
  pdfauthor={Wolf-J"urgen Beyn, Denny Otten, Jens Rottmann-Matthes},
  pdfkeywords={damped wave equation, traveling wave, freezing method, second order evolution equation, relative equilibria, stability},
  pdfcreator={latex,bibtex,dvips,ps2pdf},
  pdfproducer={LaTeX mit hyperref}
}

\newcommand{\C}{\mathbb{C}}               
\newcommand{\R}{\mathbb{R}}               
\newcommand{\Z}{\mathbb{Z}}                
\newcommand{\SO}{\mathrm{SO}}              
\newcommand{\SE}{\mathrm{SE}}              
\renewcommand{\Re}{\mathrm{Re}\,}          
\newcommand{\Q}{\mathcal{Q}}             
\newcommand{\diag}{\mathrm{diag}}          

\newcommand{\PL}{\mathcal{P}}
\newcommand{\herm}{{\mathsf{H}}}

\newcommand{\begriff}[1]{\text{#1}}

\makeatletter
\renewcommand{\@secnumfont}{\bfseries}
\makeatother
\makeatletter
  \def\section{\@startsection{section}{1}%
    \z@{.7\linespacing\@plus\linespacing}{.5\linespacing}%
    {\normalfont\LARGE\bfseries}}
\makeatother
\makeatletter
\def\@seccntformat#1{%
  \protect\textup{%
    \protect\@secnumfont
    \expandafter\protect\csname format#1\endcsname 
    \csname the#1\endcsname
    \protect\@secnumpunct
  }%
}

\titleformat*{\subsubsection}{\bfseries}

\newcommand{\sect}
{
  \setcounter{equation}{0}
  \setcounter{figure}{0}
  \section
}

\theoremstyle{plain}
\newtheorem{definition}{Definition}[section]

\theoremstyle{definition}

\newtheorem{example}[definition]{Example}
\newtheorem{proposition}[definition]{Proposition}

\allowdisplaybreaks

\begin{document}
\title[Freezing Traveling and Rotating Waves\\in Second Order Evolution Equations]{Freezing Traveling and Rotating Waves\\in Second Order Evolution Equations}
\setlength{\parindent}{0pt}
\begin{center}
\normalfont\huge\bfseries{\shorttitle}\\
\vspace*{0.25cm}
\end{center}

\vspace*{0.8cm}
\noindent
\begin{minipage}[t]{0.99\textwidth}
\begin{minipage}[t]{0.48\textwidth}
\hspace*{1.8cm} 
\textbf{Wolf-J{\"u}rgen Beyn}\footnotemark[1]${}^{,}$\footnotemark[4] \\
\hspace*{1.8cm}
\textbf{Denny Otten}\footnotemark[2]${}^{,}$\footnotemark[4] \\
\hspace*{1.8cm}
Department of Mathematics \\
\hspace*{1.8cm}
Bielefeld University \\
\hspace*{1.8cm}
33501 Bielefeld \\
\hspace*{1.8cm}
Germany
\end{minipage}
\begin{minipage}[t]{0.48\textwidth}
\hspace*{1.8cm}
\textbf{Jens Rottmann-Matthes}\footnotemark[3]${}^{,}$\footnotemark[5] \\
\hspace*{1.8cm}
Institut für Analysis \\
\hspace*{1.8cm}
Karlsruhe Institute of Technology \\
\hspace*{1.8cm}
76131 Karlsruhe \\
\hspace*{1.8cm}
Germany
\end{minipage}
\end{minipage}\\

\footnotetext[1]{e-mail: \textcolor{blue}{beyn@math.uni-bielefeld.de}, phone: \textcolor{blue}{+49 (0)521 106 4798}, \\
                                          fax: \textcolor{blue}{+49 (0)521 106 6498}, homepage: \url{http://www.math.uni-bielefeld.de/~beyn/AG\_Numerik/}.}
\footnotetext[2]{e-mail: \textcolor{blue}{dotten@math.uni-bielefeld.de}, phone: \textcolor{blue}{+49 (0)521 106 4784}, \\
                                 fax: \textcolor{blue}{+49 (0)521 106 6498}, homepage: \url{http://www.math.uni-bielefeld.de/~dotten/}.}
\footnotetext[3]{e-mail: \textcolor{blue}{jens.rottmann-matthes@kit.edu}, phone: \textcolor{blue}{+49 (0)721 608 41632}, \\
                                          fax: \textcolor{blue}{+49 (0)721 608 46530}, homepage: \url{http://www.math.kit.edu/iana2/~rottmann/}.}
\footnotetext[4]{supported by CRC 701 'Spectral Structures and Topological Methods in Mathematics',  Bielefeld University}
\footnotetext[5]{supported by CRC 1173 'Wave Phenomena: Analysis and Numerics', Karlsruhe Institute of Technology}
\vspace*{0.6cm}
\noindent
\hspace*{4.2cm}
Date: \today
\normalparindent=12pt

\vspace{0.4cm}
\begin{center}
\begin{minipage}{0.9\textwidth}
  {\small
  \textbf{Abstract.} 
  In this paper we investigate the implementation of the so-called  {\it freezing method} for second order wave equations in one and several space dimensions.
The method converts the given PDE into a partial differential algebraic
equation which is then solved numerically. The reformulation aims at 
separating the motion of a solution into a co-moving frame
and a profile which varies as little as possible. 
Numerical examples demonstrate the feasability of this approach for 
semilinear wave equations with sufficient damping. We treat the case of  a traveling 
wave in one space dimension and of a rotating wave in two space dimensions.
In addition, we investigate in arbitrary space dimensions
the point spectrum and the essential spectrum of operators obtained by linearizing
about the profile, and we indicate the consequences for the nonlinear stability of
the wave.
  }
\end{minipage}
\end{center}




%
%
\sect{Introduction}
\label{sec:1}

The topic of this paper is the numerical computation and stability of
waves occurring in nonlinear second order evolution equations with damping terms. Our main object of study is the damped wave equation
in one or several space dimensions with a nonlinearity of semilinear
type (see \eqref{equ:1.1}, \eqref{equ:1.6} below).
In the literature there are many approaches to the numerical solution of the Cauchy problem for such equations by various types of spatial and temporal discretizations. We refer, for example,
to the recent papers \cite{Brunner2013}, \cite{Glowinski2014}, \cite{Alonso2015}, \cite{Rincon2016}.
Most of the results concern finite time error estimates,
and there are a few studies of detecting blow-up solutions
or the shape of a developing solitary wave.

In our work we take a different numerical approach which
emphasizes the longtime behavior and tries to determine the
shape and speed of traveling and rotating waves from a reformulation of the
original PDE.
More specifically, we transfer the so called {\it freezing method} 
(see \cite{BeynThuemmler2004}, \cite{RowleyKevrekidisMarsden2003}, 
\cite{BeynOttenRottmannMatthes2013}) from first order to second order
evolution equations, and we investigate its relation to the stability
of the waves. Generally speaking, the method tries to separate the solution 
of a Cauchy problem into the motion of a co-moving frame and of a profile,
where the latter is required to vary as little as possible or even become 
stationary. This is achieved by transforming the original PDE into a 
partial differential algebraic equation (PDAE). The PDAE involves extra
unknowns specifying the frame, and extra constraints (so called 
{\it phase conditions}) enforcing the freezing principle for the profile. 
This methodology has been successfully applied to a wide range of PDEs 
which are of first order in time and of hyperbolic, parabolic or of mixed 
type, cf. \cite{Thuemmler2006}, \cite{Thuemmler2008}, \cite{Thuemmler2008b}, 
\cite{BeynThuemmlerSelle2008}, \cite{RottmannMatthes2012a}, \cite{RottmannMatthes2012b}, 
\cite{RottmannMatthes2012c}, \cite{BeynOttenRottmannMatthes2013}. One aim 
of the theoretical underpinning is to prove that waves which are (asymptotically) 
stable with asymptotic phase for the PDE, become stable in the classical
Lyapunov sense for the PDAE. While this has been rigorously proved
for many systems in one space dimension and confirmed numerically 
in higher space dimensions, the corresponding theory for the multi-dimensional 
case is still in its early stages, see \cite{BeynLorenz2008}, 
\cite{BeynOtten2016}, \cite{BeynOtten2016b}, \cite{Otten2015}.

In this paper we develop the freezing formulation and perform
the spectral calculations in an informal way, for the one-dimensional
as well as the multi-dimensional case. Rigorous stability results for 
the one-dimensional damped wave equation may be found in 
\cite{GallayRaugel1997}, \cite{GarrayJoly2009},
\cite{beynottenrottmann-matthes2016}.

Here we consider a nonlinear wave equation of the form
\begin{equation}
  \label{equ:1.1}
  Mu_{tt} = Au_{xx} + f(u,u_x,u_t),\,x\in\R,\,t\geqslant 0,
\end{equation}
where $u(x,t)\in \R^m, A,M \in \R^{m,m}$ and $f:\R^{3m}\to \R^m$ is 
sufficiently smooth. In addition, we assume the matrix $M$ to be nonsingular 
and $M^{-1}A$ to be positive diagonalizable, which will lead to local  
wellposedness of the Cauchy problem associated with \eqref{equ:1.1}.
Our interest is in traveling waves
\begin{equation*}
  u_{\star}(x,t) = v_{\star}(x-\mu_{\star}t),\,x\in\R,\,t\geqslant 0,
\end{equation*}
with constant limits at $\pm \infty$, i.e.
\begin{equation}
  \label{equ:1.3}
  \lim_{\xi\to\pm\infty}v_{\star}(\xi)=v_{\pm}\in\R^m,\;\lim_{\xi\to\pm\infty}v_{\star,\xi}(\xi)=0,\quad f(v_{\pm},0,0)=0.
\end{equation}

Transforming \eqref{equ:1.1} into a co-moving frame via
$u(x,t)=v(\xi,t), \xi=x-\mu_{\star}t$ leads to the system 
\begin{equation}
    \label{equ:1.4}
    Mv_{tt} = (A-\mu_{\star}^2 M)v_{\xi\xi} + 2\mu_{\star}Mv_{\xi t} + f(v,v_{\xi},v_t-\mu_{\star}v_{\xi}),\,\xi\in\R,\,t\geqslant 0.
\end{equation}
This system has $v_{\star}$ as a steady state,
\begin{equation}
  \label{equ:1.5}
  0 = (A-\mu_{\star}^2 M)v_{\star,\xi\xi} + f(v_{\star},v_{\star,\xi},-\mu_{\star}v_{\star,\xi}),\,\xi\in\R.
\end{equation}
In Section  \ref{sec:2} we work out the details of the freezing PDAE
based on the ansatz $u(x,t)=v(x-\gamma(t),t)$, $x\in \R, t\ge 0$ with the 
additional unknown function $\gamma(t), t\ge 0$. Solving this PDAE
numerically will then be demonstated for a special semilinear case,
for which damping occurs and for which the nonlinearity is of quintic 
type with $5$ zeros.
We will also discuss in Section \ref{subsec:2.2}  the spectral properties 
of the linear operator obtained by linearizing the right-hand side of 
\eqref{equ:1.4} about the profile $v_{\star}$. First, there is the eigenvalue 
zero due to shift equivariance, and then we analyze the dispersion curves 
which are part of the operator's essential spectrum. If there is sufficient
damping in the system (depending on the derivative $D_3f$), one can expect 
the whole nonzero spectrum to lie strictly to the left of the imaginary axis. 
We refer to \cite{beynottenrottmann-matthes2016} for a rigorous proof of  
nonlinear stability in such a situation, both stability of the wave with 
asymptotic phase for equation \eqref{equ:1.4} and Lyapunov stability of 
the wave and its speed for the freezing equation.
 
The subsequent section is devoted to study corresponding problems
for multi-dimensional wave equations
\begin{equation} 
  \label{equ:1.6}
  M u_{tt} + Bu_t= A \Delta u +f(u),\,x\in \R^d,\,t\geqslant 0,
\end{equation}
where the matrices $A,M$ are as above, the damping matrix $B \in \R^{m,m}$ is
given  and $f:\R^m\to \R^m$ is again
sufficiently smooth. We look for rotating waves of the form
\begin{equation*} 
  u_{\star}(x,t) = v_{\star}(e^{-tS_{\star}}(x-x_{\star})),\,x\in\R^d,\,t\geqslant 0,
\end{equation*}
where $x_{\star}\in \R^d$ denotes the center of rotation, $S_{\star}\in \R^{d,d}$
is a skew-symmetric matrix, and $v_{\star}:\R^d \to \R^m$ describes
the profile. Transforming \eqref{equ:1.6} into a co-rotating frame
via $u(x,t)=v(e^{-tS_{\star}}(x-x_{\star}),t)$ now leads to the equation
\begin{equation}
\label{equ:1.8}
  \begin{aligned}
  Mv_{tt}+Bv_t =& A\triangle v - Mv_{\xi\xi}(S_{\star}\xi)^2 +2Mv_{\xi t}S_{\star}\xi  - Mv_{\xi}S_{\star}^2\xi 
                + Bv_{\xi}S_{\star}\xi + f(v),\,\xi\in\R^d,\,t\geqslant 0,
  \end{aligned}
\end{equation}
where our notation for derivatives uses multilinear calculus, e.g.
\begin{align}
  \label{equ:1.8a}
  (v_{\xi\xi}h_1 h_2)_i = \sum_{j=1}^{d}\sum_{k=1}^{d}v_{i,\xi_j\xi_k}(h_1)_j(h_2)_k,\quad (\triangle v)_i = \sum_{j=1}^{d}v_{i,\xi_j\xi_j} = \sum_{j=1}^{d}v_{i,\xi\xi}(e^j)^2. 
\end{align}
The profile $v_{\star}$ of the wave is then a steady state solution of \eqref{equ:1.8}, i.e.
\begin{equation}
  \label{equ:1.9}
    0 = A\triangle v_{\star} - Mv_{\star,\xi\xi}(S_{\star}\xi)^2- 
Mv_{\star,\xi}S_{\star}^2\xi +Bv_{\star,\xi}S_{\star}\xi  + f(v_{\star}),\,\xi\in\R^d.
\end{equation}
As is known from first oder in time PDEs, there are several eigenvalues of
the linearized operator on the imaginary axis caused by the Euclidean symmmetry, 
see e.g. \cite{Metafune2001}, \cite{MetafunePallaraPriola2002}, 
\cite{FiedlerScheel2003}, \cite{BeynLorenz2008}, \cite{Otten2014}. 
The computations become more involved for the wave equation \eqref{equ:1.8}, 
but we will show that the eigenvalues on the imaginary axis are
the same as in the parabolic case. 
 Further, determining the dispersion relation, and thus curves 
in the essential spectrum, now amounts to solving a parameterized quadratic 
eigenvalue problem which in general can only be solved numerically. Finally, 
we present a numerical example of a rotating wave for the cubic-quintic 
Ginzburg-Landau equation. The performance of the freezing method will be 
demonstrated, and we investigate the numerical eigenvalues approximating 
the point spectrum on (and close to) the imaginary axis as well as the 
essential spectrum in the left half-plane.

%
%
\sect{Traveling waves in one space dimension}
\label{sec:2}

\subsection{Freezing traveling waves.}
\label{subsec:2.1}

Consider the Cauchy problem associated with \eqref{equ:1.1}
\begin{subequations} 
  \label{equ:2.1}
  \begin{align}
    & Mu_{tt} = Au_{xx} + f(u,u_x,u_t),           &&\,x\in\R,\,t\geqslant 0, \label{equ:2.1a} \\
    & u(\cdot,0) = u_0,\quad u_t(\cdot,0) = v_0,  &&\,x\in\R,\,t=0, \label{equ:2.1b}
  \end{align}
\end{subequations}
for some initial data $u_0,v_0:\R\rightarrow\R^m$ and some nonlinearity $f\in C^3(\R^{3m},\R)$. 
Introducing new unknowns $\gamma(t)\in\R$ and $v(\xi,t)\in\R^m$ via the \begriff{freezing ansatz for traveling waves}
\begin{equation}
  \begin{aligned}
  \label{equ:2.2}
    u(x,t) & = v(\xi,t),\quad\xi:=x-\gamma(t),\,x\in\R,\,t\geqslant 0,
  \end{aligned}
\end{equation}
and inserting \eqref{equ:2.2} into \eqref{equ:2.1a} by taking
\begin{align}
  \label{equ:2.3}
  u_t = -\gamma_t v_{\xi}+v_t,\quad u_{tt}=-\gamma_{tt}v_{\xi}+\gamma_t^2v_{\xi\xi}-2\gamma_t v_{\xi t} + v_{tt}
\end{align}
into account, we obtain the equation
\begin{equation}
  \label{equ:2.4}
  Mv_{tt} = (A-\gamma_t^2 M)v_{\xi\xi} + 2\gamma_t M v_{\xi t} + \gamma_{tt}M v_{\xi} + f(v,v_{\xi},v_t-\gamma_t v_{\xi}),\;\xi\in\R,\,t\geqslant 0.
\end{equation}
Now it is convenient to introduce time-dependent functions $\mu_1(t)\in\R$ and $\mu_2(t)\in\R$ via
\begin{equation*}
  \mu_1(t):=\gamma_t(t),\quad \mu_2(t):=\mu_{1,t}(t)=\gamma_{tt}(t)
\end{equation*}
which allows us to transfer \eqref{equ:2.4} into a coupled PDE/ODE-system
\begin{subequations} 
  \label{equ:2.5}
  \begin{align}
    &Mv_{tt} = (A-\mu_1^2 M)v_{\xi\xi} + 2\mu_1 M v_{\xi t} + \mu_2 M v_{\xi} + f(v,v_{\xi},v_t-\mu_1 v_{\xi}), &&\xi\in\R,\,t\geqslant 0,  \label{equ:2.5a}\\
    &\mu_{1,t} = \mu_2, &&t\geqslant 0,\label{equ:2.5b}\\
    &\gamma_t = \mu_1, &&t\geqslant 0.\label{equ:2.5c}
  \end{align}
\end{subequations}
The quantity $\gamma(t)$ denotes the \begriff{position}, $\mu_1(t)$ the velocity and $\mu_2(t)$ the acceleration of the profile $v(\xi,t)$ at time $t$. 
We next specify initial data for the system \eqref{equ:2.5} as follows,
\begin{equation}
  \label{equ:2.6}
  v(\cdot,0) = u_0,\quad v_t(\cdot,0)=v_0+\mu_1^0 u_{0,\xi},\quad \mu_1(0)=\mu_1^0,\quad \gamma(0)=0
\end{equation}
Note that if we require $\gamma(0)=0$ and $\mu_1(0)=\mu_1^0$, then the first equation in \eqref{equ:2.6} follows from \eqref{equ:2.2} and \eqref{equ:2.1b}, 
while the second equation in \eqref{equ:2.6} follows from \eqref{equ:2.3}, \eqref{equ:2.1b} and \eqref{equ:2.5c}. Suitable values for $\mu_1^0$ depend on 
the choice of phase condition to be discussed next.

We compensate the extra variable $\mu_2$ in the system \eqref{equ:2.5} by imposing an additional scalar algebraic constraint, also known as a phase condition, 
of the general form
\begin{equation}
  \label{equ:2.7}
  \psi(v,v_t,\mu_1,\mu_2) = 0,\;t\geqslant 0.
\end{equation}
Two possible choices are the \begriff{fixed phase condition} $\psi_{\mathrm{fix}}$ and the \begriff{orthogonal phase condition} $\psi_{\mathrm{orth}}$ given by
\begin{align}
  \psi_{\mathrm{fix}}(v) = \langle v-\hat{v},\hat{v}_{\xi}\rangle_{L^2},\;t\geqslant 0, \label{equ:2.8} \\
  \psi_{\mathrm{orth}}(v_t) = \langle v_t,v_{\xi}\rangle_{L^2},\;t\geqslant 0. \nonumber
\end{align}
These two types and their derivation are discussed in \cite{beynottenrottmann-matthes2016}. The function $\hat{v}:\R\rightarrow\R^m$ denotes a 
time-independent and sufficiently smooth template (or reference) function, e.g. $\hat{v}=u_0$. Suitable values for $\mu_1(0)=\mu_1^0$ can be derived 
from requiring consistent initial values for the PDAE. For example, consider \eqref{equ:2.8} and take the time derivative at $t=0$. Together with \eqref{equ:2.6} 
this leads to $0=\langle v_t(\cdot,0),\hat{v}_{\xi}\rangle_{L^2}=\langle v_0,\hat{v}_{\xi}\rangle_{L^2}+\mu_1^0\langle u_{0,\xi},\hat{v}_{\xi}\rangle_{L^2}$. 
If $\langle u_{0,\xi},\hat{v}_{\xi}\rangle_{L^2}\neq 0$ this determines a unique value for $\mu_1^0$.

Let us summarize the set of equations obtained by the freezing method of the original Cauchy problem \eqref{equ:2.1}. Combining the differential equations 
\eqref{equ:2.5}, the initial data \eqref{equ:2.6} and the phase condition \eqref{equ:2.7}, we arrive at the following partial differential algebraic evolution 
equation (short: PDAE) to be solved numerically:
\begin{subequations}
  \label{equ:2.10}
  \begin{align}
     \label{equ:2.10a}
     &\begin{aligned}
      Mv_{tt} &= (A-\mu_1^2 M)v_{\xi\xi} + 2\mu_1 M v_{\xi,t} + \mu_2 M v_{\xi} + f(v,v_{\xi},v_t-\mu_1 v_{\xi}),\\
      \mu_{1,t} &= \mu_2, \quad \gamma_t=\mu_1,
    \end{aligned}&t\geqslant 0,\\
    & 0 = \psi(v,v_t,\mu_1,\mu_2), &t\geqslant 0,\label{equ:2.10b}\\
    &\begin{aligned}
      v(\cdot,0) &= u_0,\quad v_t(\cdot,0) = v_0+\mu_1^0 u_{0,\xi}, \quad \mu_1(0) = \mu_1^0, \quad \gamma(0) = 0.
    \end{aligned}\label{equ:2.10c}
  \end{align}
\end{subequations}
The system \eqref{equ:2.10} depends on the choice of phase condition $\psi$ and is to be solved for $(v,\mu_1,\mu_2,\gamma)$ 
with given initial data $(u_0,v_0,\mu_1^0)$. It consists of a PDE for $v$ that is coupled to two ODEs for $\mu_1$ and $\gamma$ 
\eqref{equ:2.10a} and an algebraic constraint \eqref{equ:2.10b} which closes the system. A consistent initial value $\mu_1^0$ 
for $\mu_1$ is computed from the phase condition and the initial data. Further initialization of the algebraic variable $\mu_2$ 
is usually not needed for a PDAE-solver but can be provided if necessary (see \cite{beynottenrottmann-matthes2016}).

The ODE for $\gamma$ is called the \begriff{reconstruction equation} in \cite{RowleyKevrekidisMarsden2003}. It decouples from 
the other equations in \eqref{equ:2.10} and can be solved in a postprocessing step. The ODE for $\mu_1$ is the new feature of 
the PDAE for second order systems when compared to the first order parabolic and hyperbolic equations in 
\cite{BeynThuemmler2004,RottmannMatthes2010,BeynOttenRottmannMatthes2013}.

Finally, note that $(v,\mu_1,\mu_2)=(v_{\star},\mu_{\star},0)$ satisfies
\begin{align*}
  0 & = (A-\mu_{\star}^2 M)v_{\star,\xi\xi} + \mu_{\star} M v_{\star,\xi} + f(v_{\star},v_{\star,\xi},-\mu_{\star}v_{\star,\xi}),\;\xi\in\R,\\
  0 & = \mu_2, \\
  0 & = \psi(v_{\star},0,\mu_{\star},0),
\end{align*}
and hence is a stationary solution of \eqref{equ:2.10a},\eqref{equ:2.10b}. Here we assume that $v_{\star},\mu_{\star}$ have been selected 
to satisfy the phase condition. Obviously, in this case we have $\gamma(t)=\mu_{\star}t$. 
For a stable traveling wave we expect that solutions $(v,\mu_1,\mu_2,\gamma)$ of \eqref{equ:2.10} show the limiting behavior
\begin{align*}
  v(t)\rightarrow v_{\star},\quad \mu_1(t)\rightarrow\mu_{\star},\quad \mu_2(t)\rightarrow 0\quad\text{as}\quad t\to\infty,
\end{align*}
provided the initial data are close to their limiting values. 

\begin{example}[Freezing quintic Nagumo wave equation]\label{exa:1}
  Consider the quintic Nagumo wave equation, 
  \begin{equation}
    \label{equ:2.22}
    \varepsilon u_{tt}  = A u_{xx} + f(u,u_x,u_t),\; x\in \R,\, t\geqslant 0,
  \end{equation}
  with $u=u(x,t)\in\R$, $\varepsilon>0$,  $0<\alpha_1<\alpha_2<\alpha_3<1$,
and the nonlinear term
  \begin{equation} \label{equ:2.22a}
    f:\R^3\rightarrow\R,\quad f(u,u_x,u_t)=-u_t+u(1-u)\prod_{j=1}^{3}(u-\alpha_j).
  \end{equation}
  
  \begin{figure}[ht]
    \centering
    \subfigure[]{\includegraphics[height=3.8cm] {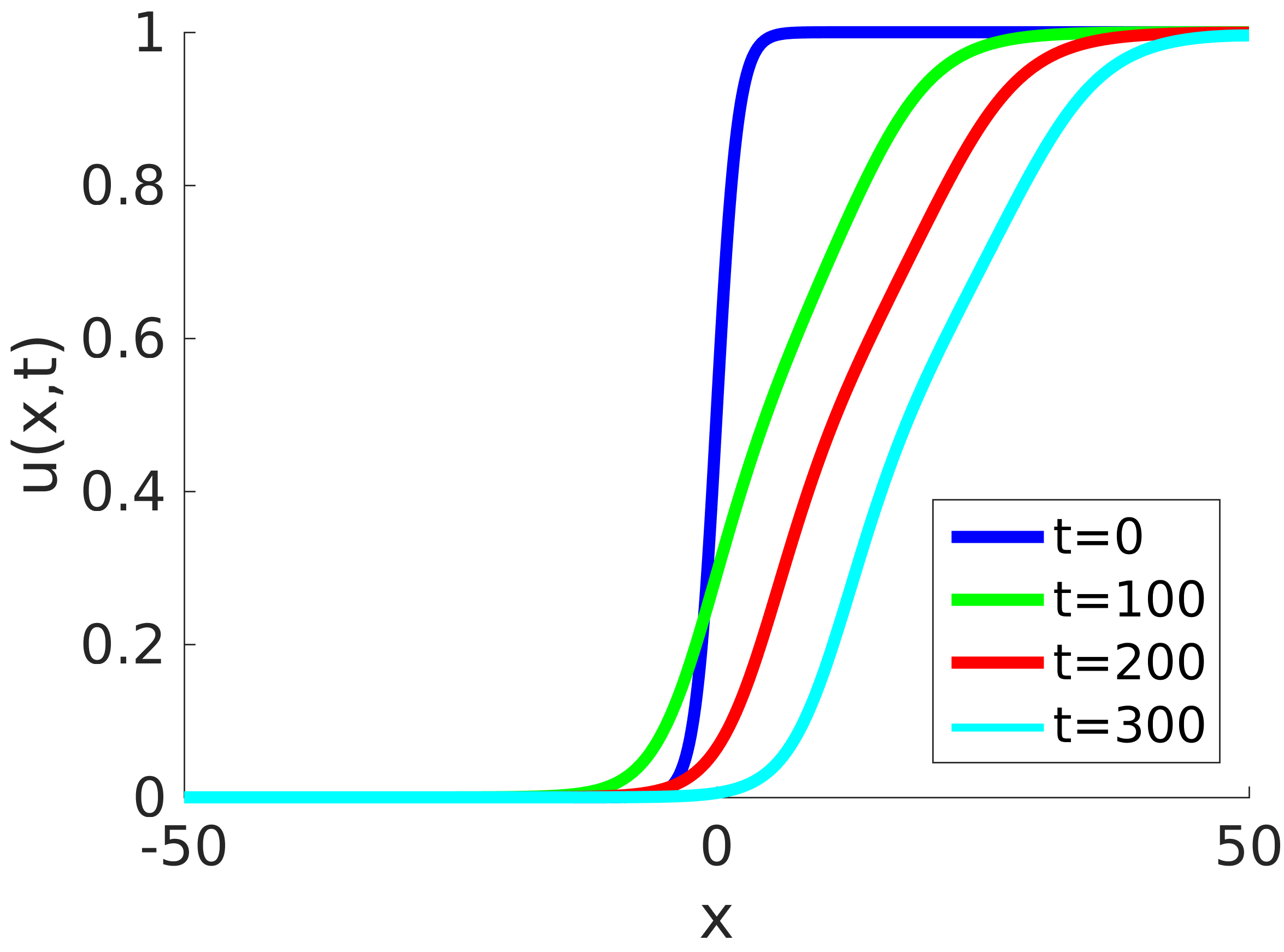}\label{fig:2.1a}}
    \subfigure[]{\includegraphics[height=3.8cm] {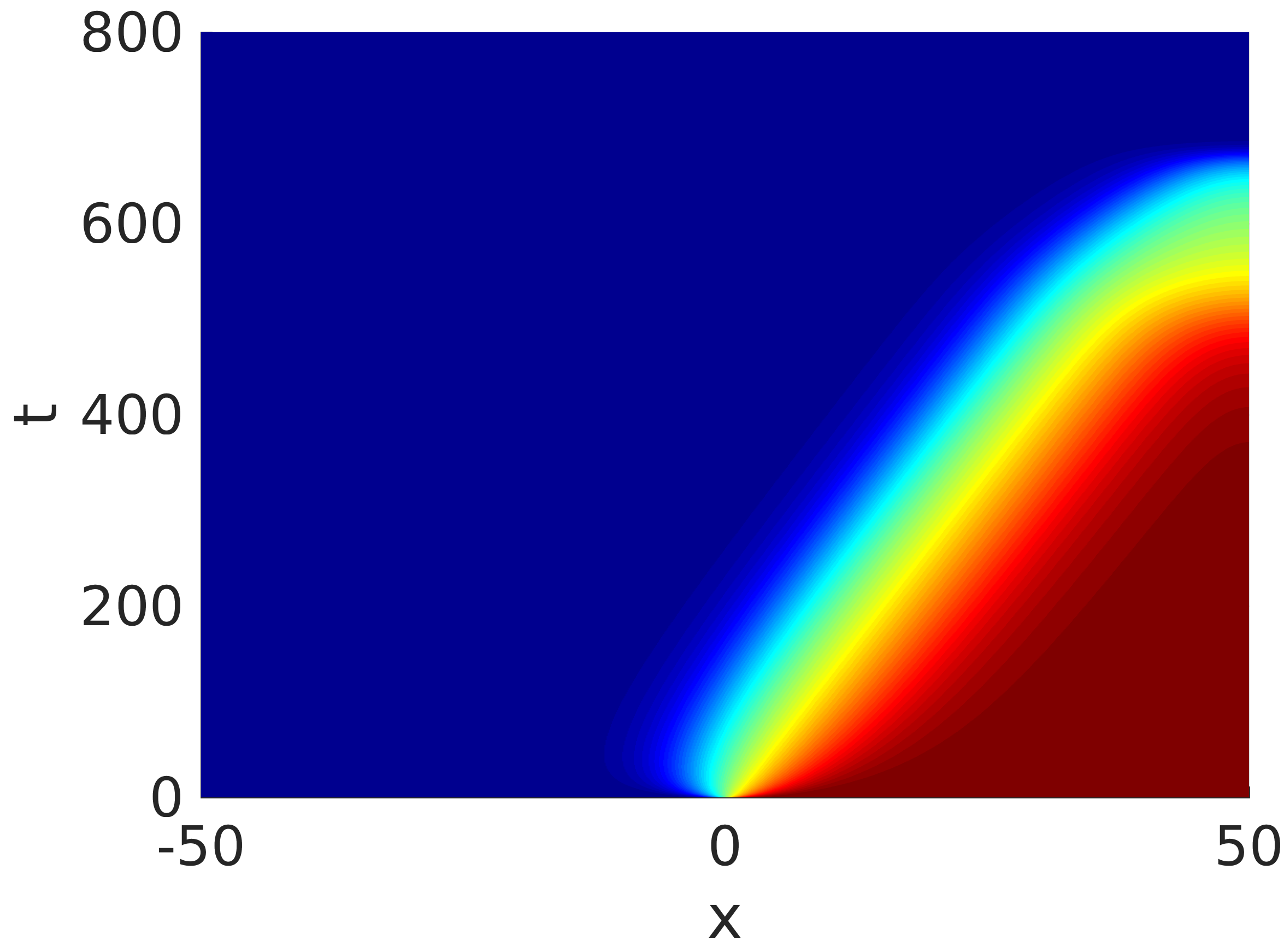} \label{fig:2.1b}}
    \caption{Traveling front of quintic Nagumo wave equation \eqref{equ:2.22} at different time instances (a) and its time evolution (b) for parameters from \eqref{equ:2.23}.}
    \label{fig:2.1}
  \end{figure}
For the parameter values
  \begin{equation}
    \label{equ:2.23}
    M=\varepsilon=\frac{1}{2},\quad A=1,\quad\alpha_1=\frac{2}{5},\quad \alpha_2=\frac{1}{2},\quad \alpha_3=\frac{17}{20},
  \end{equation}
  equation \eqref{equ:2.22} admits a traveling front solution connecting the asymptotic states $v_-=0$ and $v_+=1$.
 
 Figure \ref{fig:2.1} shows a numerical simulation of the solution $u$ of \eqref{equ:2.22} on the spatial domain $(-50,50)$ with homogeneous 
  Neumann boundary conditions, with initial data
  \begin{align} 
    \label{equ:2.24} 
    u_0(x)=\tfrac{1}{2}\left(1+\tanh\left(\tfrac{x}{2}\right)\right),\quad v_0(x)=0
  \end{align}
  and parameters taken from \eqref{equ:2.23}. For the space discretization we use continuous piecewise linear finite 
  elements with spatial stepsize $\triangle x=0.1$. For the time discretization we use the BDF method of order $2$ 
  with absolute tolerance $\mathrm{atol}=10^{-3}$, relative tolerance $\mathrm{rtol}=10^{-2}$, temporal stepsize $\triangle t=0.2$
  and final time $T=800$. Computations are performed with the help of the software COMSOL 5.2.
   
  Let us now consider the frozen quintic Nagumo wave equation resulting from \eqref{equ:2.10}
  \begin{subequations}    
  \label{equ:2.25}
  \begin{align}
     \label{equ:2.25a}
     &\begin{aligned}
      \varepsilon v_{tt} + v_t &= (1-\mu_1^2 \varepsilon)v_{\xi\xi} +
      2\mu_1 \varepsilon v_{\xi,t} +
      (\mu_2 \varepsilon + \mu_1)v_{\xi} + \tilde{f}(v),\\
      \mu_{1,t} &= \mu_2, \quad \gamma_t=\mu_1,
    \end{aligned}&t\geqslant 0,\\
    & 0 = \bigl\langle v_t(\cdot,t),\hat{v}_\xi\bigr\rangle_{L^2(\R,\R)}, &t\geqslant 0,\label{equ:2.25b}\\
    &\begin{aligned}
      v(\cdot,0) &= u_0,\quad v_t(\cdot,0) = v_0+\mu_1^0 u_{0,\xi}, \quad
      \mu_1(0) = \mu_1^0, \quad
      \gamma(0) = 0.
    \end{aligned}\label{equ:2.25c}
  \end{align}
  \end{subequations}
  
  \begin{figure}[ht]
  \centering
  \subfigure[]{\includegraphics[height=3.8cm] {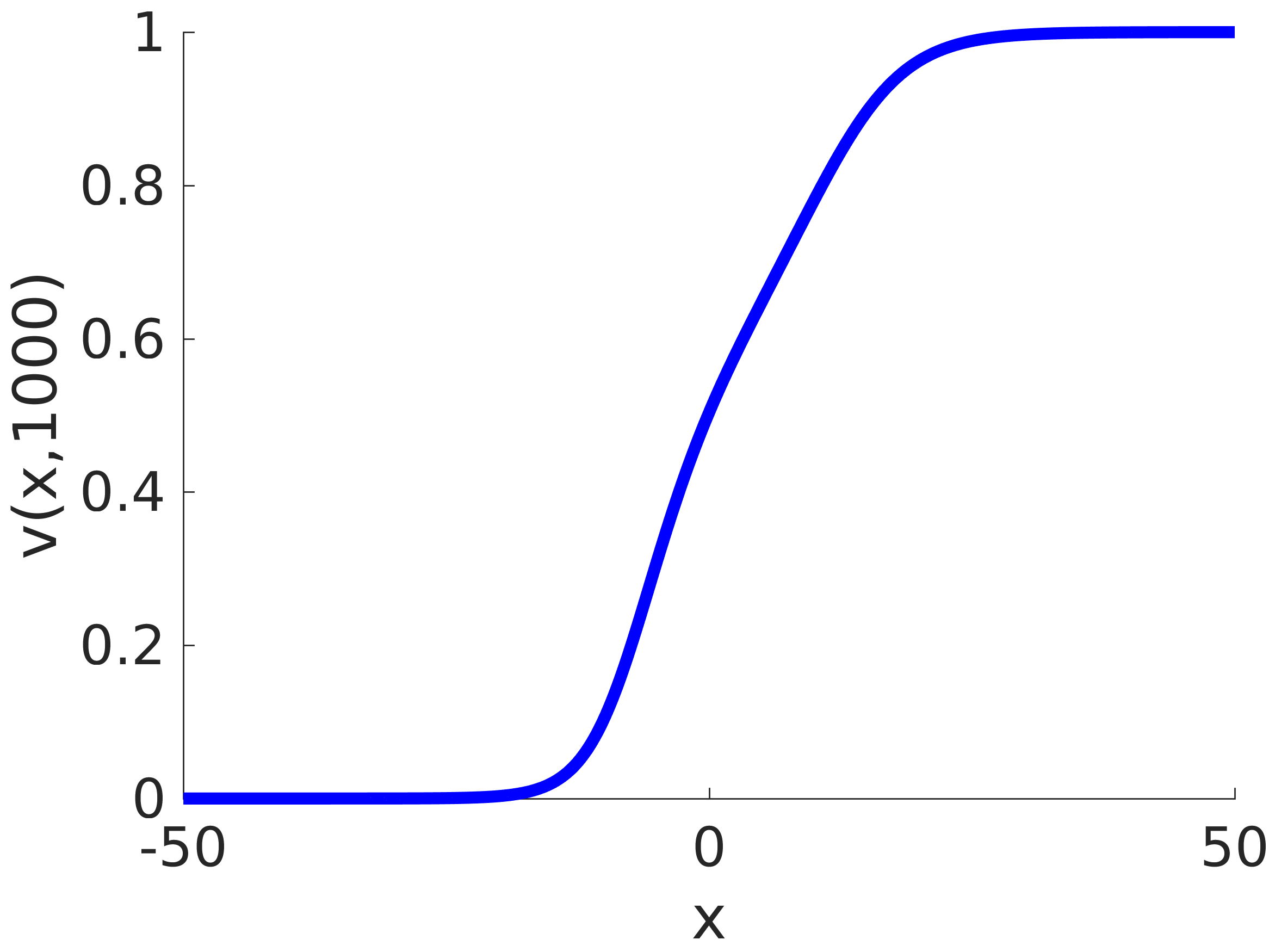}\label{fig:2.2a}}
  \subfigure[]{\includegraphics[height=3.8cm] {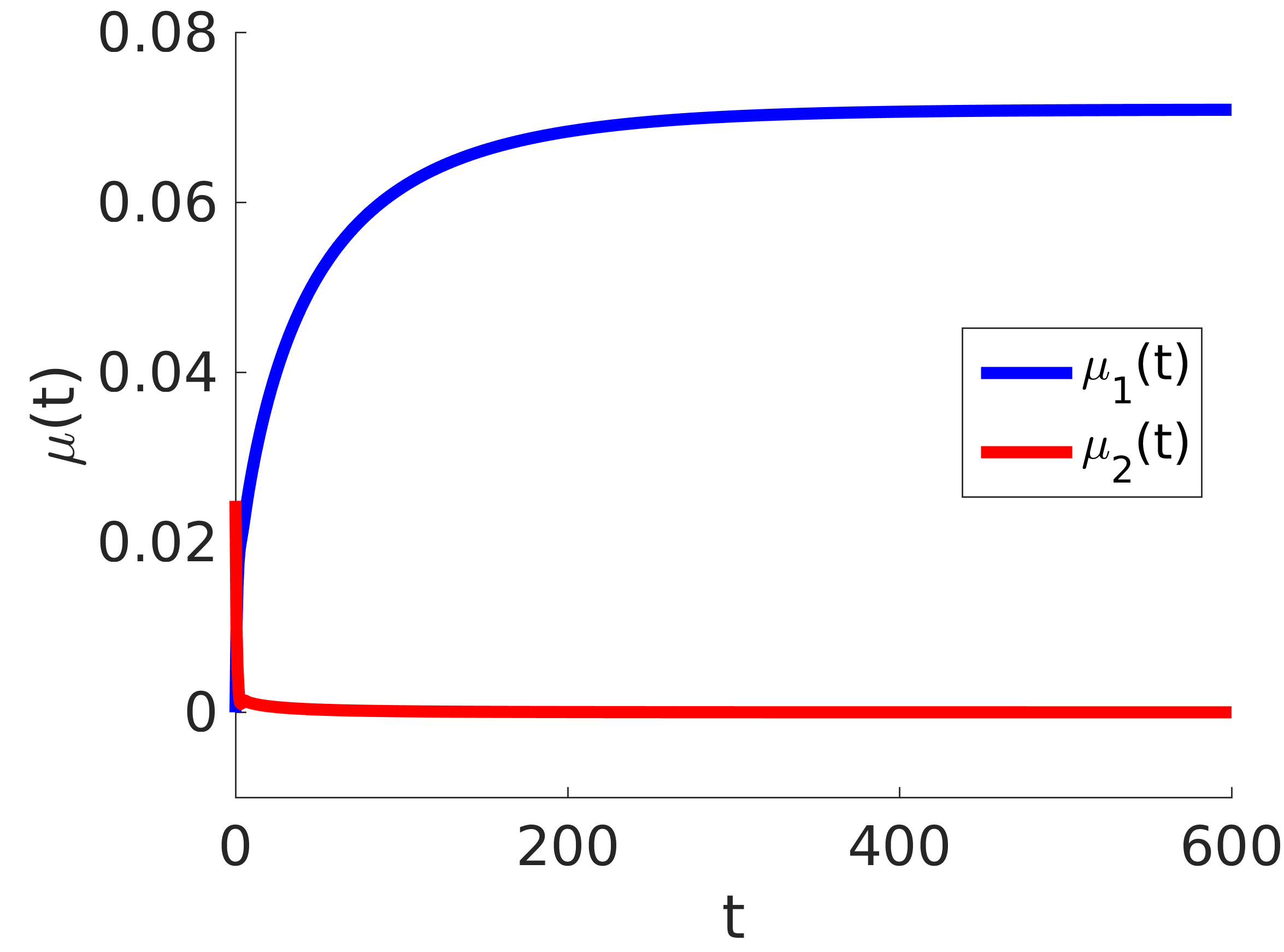} \label{fig:2.2b}}
  \subfigure[]{\includegraphics[height=3.8cm] {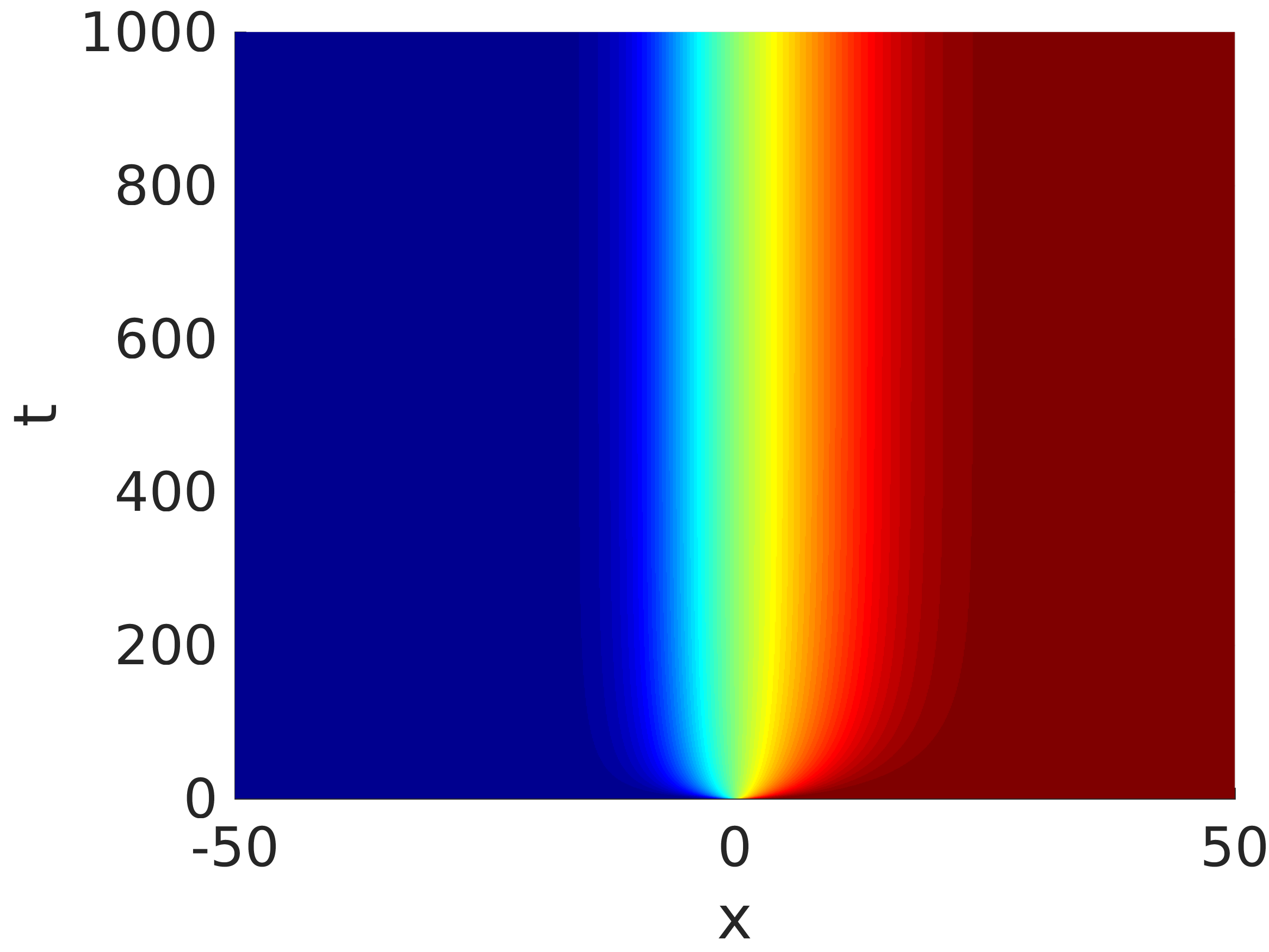}\label{fig:2.2c}}
  \caption{Solution of the frozen quintic Nagumo wave equation \eqref{equ:2.25}: approximation of profile $v(x,1000)$ (a) and time evolutions of velocity $\mu_1$ 
  and acceleration $\mu_2$ (b) and of the profile $v$ (c) for parameters from \eqref{equ:2.23}.}
  \label{fig:2.2}
  \end{figure}

  Figure \ref{fig:2.2} shows the solution $(v,\mu_1,\mu_2,\gamma)$ of \eqref{equ:2.25} on the spatial domain $(-50,50)$ with homogeneous Neumann 
  boundary conditions, initial data $u_0$, $v_0$ from \eqref{equ:2.24}, and reference function $\hat{v}=u_0$. For the computation we used the 
  fixed phase condition $\psi_{\mathrm{fix}}(v)$ from \eqref{equ:2.8} with consistent intial data $\mu_1^0=0$, see above. 
  The spatial discretization data are taken as in the nonfrozen case. For the time discretization we used the BDF method of order $2$ 
  with absolute tolerance $\mathrm{atol}=10^{-3}$, relative tolerance $\mathrm{rtol}=10^{-2}$, temporal stepsize $\triangle t=0.6$
  and final time $T=3000$. The diagrams show that after a very short transition phase the profile becomes stationary, the acceleration $\mu_2$ 
  converges to zero, and the speed $\mu_1$ approaches an asymptotic value $\mu_{\star}^{\mathrm{num}}= 0.0709$.
\end{example}

\subsection{Spectra of traveling waves.}
\label{subsec:2.2}

Consider the linearized equation 
\begin{equation}
  \label{equ:2.11}
  M v_{tt} - (A-\mu_{\star}^2 M)v_{\xi\xi} - 2\mu_{\star}Mv_{\xi t} - (D_2 f_{\star} - \mu_{\star} D_3 f_{\star})v_{\xi} - D_3 f_{\star} v_t - D_1 f_{\star} v = 0
\end{equation}
which is obtained from the co-moving frame \eqref{equ:1.4} linearized at the profile $v_{\star}$. In \eqref{equ:2.11} we use the short form 
$D_j f_{\star}=D_j f(v_{\star},v_{\star,\xi},-\mu_{\star}v_{\star,\xi})$. 
Looking for solutions of the form $v(\xi,t)=e^{\lambda t}w(\xi)$ to \eqref{equ:2.11} yields the quadratic eigenvalue problem
\begin{equation}
  \label{equ:2.12}
  \PL(\lambda)w = \left(\lambda^2 P_2 + \lambda P_1 + P_0\right)w = 0,\,\xi\in\R
\end{equation}
with differential operators $P_j$ defined by
\begin{align*}
  P_2 = M,\quad 
  P_1 = -2\mu_{\star}M\partial_\xi - D_3 f_{\star},\quad
  P_0 = -(A-\mu_{\star}^2 M)\partial^2_\xi - (D_2 f_{\star} - \mu_{\star}D_3 f_{\star})\partial_\xi - D_1 f_{\star}.
\end{align*}
We are interested in solutions $(\lambda,w)$ of \eqref{equ:2.12} which are
candidates for eigenvalues $\lambda\in\C$ and eigenfunctions $w:\R\to\C^m$
in suitable function spaces. 
In fact, it is usually imposssible to determine  the spectrum $\sigma(\PL)$ analytically, but one is able to analyze certain 
subsets. Let us first calculate the symmetry set $\sigma_{\mathrm{sym}}(\PL)$, which belongs to the point spectrum $\sigma_{\mathrm{pt}}(\PL)$ 
and is affected by the underlying group symmetries. Then, we calculate the dispersion set $\sigma_{\mathrm{disp}}(\PL)$, which belongs to the 
essential spectrum $\sigma_{\mathrm{ess}}(\PL)$ and is affected by the far-field behavior of the wave. Let us first derive the symmetry set of $\PL$. 
This is a simple task for traveling waves but becomes more involved when analyzing the symmetry set for rotating waves (see Section \ref{subsubsec:3.2.1}).

\subsubsection{Point Spectrum and symmetry set.} 
Applying $\partial_\xi$ to the traveling wave equation \eqref{equ:1.5} yields $P_0 v_{\star,\xi}=0$ which proves the following result.

\begin{proposition}[Point spectrum of traveling waves]\label{prop:2.1}
  Let $f\in C^1(\R^{3m},\R^m)$ and let $v_{\star}\in C^3(\R,\R^m)$ be a nontrivial classical solution of \eqref{equ:1.5} for some $\mu_{\star}\in\R$. 
  Then, $w=v_{\star,\xi}$ and $\lambda=0$ is a classical solution of the eigenvalue problem \eqref{equ:2.12}. In particular, the symmetry set
  \begin{align*}
    \sigma_{\mathrm{sym}}(\PL)=\{0\}
  \end{align*}
  belongs to the point spectrum $\sigma_{\mathrm{pt}}(\PL)$ of $\PL$.
\end{proposition}
Of course, a rigorous statement of this kind requires to specify the
function spaces involved, e.g. $L^2(\R,\R^m)$ or $H^1(\R,\R^m)$, see
\cite{GallayRaugel1997}, \cite{GarrayJoly2009}, \cite{beynottenrottmann-matthes2016}.

\subsubsection{Essential Spectrum and dispersion set.}
\label{subsubsec:2.2.2}
\begin{enumerate}[label=\bf{\arabic*.},leftmargin=*]
  \item \textbf{The far-field operator.} 
It is a  well known fact that the essential spectrum is affected by the limiting equation obtained from \eqref{equ:2.12} as $\xi\to\pm\infty$. 
Therefore, we let formally $\xi\to\pm\infty$  in \eqref{equ:2.12}
and obtain
\begin{equation}
  \label{equ:2.15}
  \left(\lambda^2 P_2 + \lambda P_1^{\pm} + P_0^{\pm}\right)w=0,\;\xi\in\R.
\end{equation}
with the constant coefficient operators
\begin{align*}
  P_2 = M,\quad 
  P_1^{\pm} = -2\mu_{\star}M\partial_\xi - D_3 f_{\pm},\quad
  P_0^{\pm} = -(A-\mu_{\star}^2 M)\partial^2_\xi - (D_2 f_{\pm} - \mu_{\star}D_3 f_{\pm})\partial_\xi - D_1 f_{\pm},
\end{align*}
where $v_{\pm}$ are from \eqref{equ:1.3} and
$D_j f_{\pm}=D_j f(v_{\pm},0,0)$.
We may then write equation \eqref{equ:2.12} as
\begin{equation*}
  \left(\lambda^2 P_2 + \lambda (P_1^{\pm}+Q_1^{\pm}(\xi)) + (P_0^{\pm} + Q_2^{\pm}(\xi)\partial_\xi + Q_3^{\pm}(\xi))\right)w=0,\;\xi\in\R
\end{equation*}
 with the perturbation operators defined by
\begin{equation*}
  Q_1^{\pm}(\xi)=D_3 f_{\pm}-D_3f_{\star},\;\;
  Q_2^{\pm}(\xi)=D_2 f_{\pm} -D_2f_{\star} + \mu_{\star}(D_3 f_{\star}- D_3 f_{\pm}),\;\;
  Q_3^{\pm}(\xi)=D_1 f_{\pm}-D_1f_{\star},
\end{equation*} 
Note that $v_{\star}(\xi)\to v_{\pm}$ implies $Q_j^{\pm}(\xi)\to 0$ as $\xi\to\pm\infty$ for $j=1,2,3$.
\item \textbf{Spatial Fourier transform.} 
For $\omega\in\R$, $z\in\C^m$, $|z|=1$ we apply the spatial Fourier transform $w(\xi)=e^{i\omega\xi}z$ to equation \eqref{equ:2.15} which leads 
to the $m$-dimensional quadratic eigenvalue problem
\begin{equation}
  \label{equ:2.16}
  \left(\lambda^2 A_2 + \lambda A_1^{\pm}(\omega) + A_0^{\pm}(\omega)\right)z = 0
\end{equation}
with matrices $A_2\in\R^{m,m}$ and $A_1^{\pm}(\omega),A_0^{\pm}(\omega)\in\C^{m,m}$ given by
\begin{equation}
  \label{equ:2.16a}
  A_2 = M,\;
  A_1^{\pm}(\omega) = -2i\omega\mu_{\star}M - D_3 f_{\pm},\;
  A_0^{\pm}(\omega) = \omega^2(A-\mu_{\star}^2 M) - i\omega(D_2 f_{\pm} - \mu_{\star} D_3 f_{\pm}) - D_1 f_{\pm}.
\end{equation}
\item \textbf{Dispersion relation and dispersion set.} The dispersion relation for traveling waves of second order evolution equations states the following: 
Every $\lambda\in\C$ satisfying
\begin{equation}
  \label{equ:2.17}
  \det\left(\lambda^2 A_2 + \lambda A_1^{\pm}(\omega) + A_0^{\pm}(\omega)\right)=0
\end{equation}
for some $\omega\in\R$ belongs to the essential spectrum of $\PL$, i.e. $\lambda\in\sigma_{\mathrm{ess}}(\PL)$. Solving \eqref{equ:2.17} is equivalent 
to finding all zeros of a polynomial of degree $2m$. Note that the limiting case $M=0$ in \eqref{equ:2.17} leads to the dispersion relation for traveling 
waves of first order evolution equations, which is well-known in the literature, see \cite{Sandstede2002}.
\end{enumerate}

\begin{proposition}[Essential spectrum of traveling waves]\label{prop:2.2}
  Let $f\in C^1(\R^{3m},\R^m)$ with $f(v_{\pm},0,0)=0$ for some $v_{\pm}\in\R^m$. Let $v_{\star}\in C^2(\R,\R^m)$, $\mu_{\star}\in\R$ 
  be a nontrivial classical solution of \eqref{equ:1.5} satisfying $v_{\star}(\xi)\to v_{\pm}$ as $\xi\to\pm\infty$. Then, the dispersion set
  \begin{align*}
    \sigma_{\mathrm{disp}}(\PL) = \{\lambda\in\C : \text{$\lambda$ satisfies \eqref{equ:2.17} for some $\omega\in\R$, and $+$ or $-$}\}
  \end{align*}
  belongs to the essential spectrum $\sigma_{\mathrm{ess}}(\PL)$ of $\PL$.
\end{proposition}

\begin{example}[Spectrum of quintic Nagumo wave equation]\label{exa:2}
  As shown in Example \ref{exa:1} the quintic Nagumo wave equation \eqref{equ:2.22} with coefficients and parameters \eqref{equ:2.23}
  has a traveling front solution $u_{\star}(x,t)=v_{\star}(x-\mu_{\star}t)$ with velocity $\mu_{\star}\approx 0.0709$, whose
  profile $v_{\star}$ connects the asymptotic states $v_{-}=0$ and $v_{+}=1$ according to \eqref{equ:1.3}.\\
  We solve numerically the eigenvalue problem for the quintic Nagumo wave equation
  \begin{align}
  \label{equ:2.18a}
  \left(\lambda^2\varepsilon + \lambda\left(-2\mu_{\star}\varepsilon\partial_\xi - D_3 f_{\star}\right) + \left(-(1-\mu_{\star}^2 \varepsilon)\partial^2_\xi   
  - (D_2 f_{\star} - \mu_{\star}D_3 f_{\star})\partial_\xi - D_1 f_{\star}\right)\right)w = 0.
  \end{align}
  Both approximations of the profile $v_{\star}$ and the velocity $\mu_{\star}$ in \eqref{equ:2.18a} are chosen from the solution of \eqref{equ:2.25} at time $t=3000$ in Example \ref{exa:1}.   
  Due to Proposition \ref{prop:2.1} we expect $\lambda=0$ to be an isolated eigenvalue belonging to the point spectrum. Let us next discuss 
  the dispersion set from Proposition \ref{prop:2.2}. The quintic Nagumo 
nonlinearity  \eqref{equ:2.22a} satisfies
  \begin{equation*}
    f_{\pm}=0,\quad D_3f_{\pm}=-1,\quad D_2f_{\pm}=0,\quad D_1f_{-}=-\alpha_1\alpha_2\alpha_3,\quad D_1f_{+}=-\prod_{j=1}^{3}(1-\alpha_j).
  \end{equation*}
  The matrices $A_2$, $A_1^{\pm}(\omega)$, $A_0^{\pm}(\omega)$ from \eqref{equ:2.16a} of the quadratic problem \eqref{equ:2.16} are given by
  \begin{equation*}
    A_2=\varepsilon,\quad A_1^{\pm}(\omega)=-2i\omega\mu_{\star}\varepsilon+1,\quad A_0^{\pm}(\omega)=\omega^2(1-\mu_{\star}^2\varepsilon)-i\omega\mu_{\star}-D_1 f_{\pm}.
  \end{equation*}
  The dispersion relation \eqref{equ:2.17} for the quintic Nagumo front  states that every $\lambda\in\C$ satisfying
  \begin{equation}
    \label{equ:2.20}
    \lambda^2\varepsilon + \lambda(-2i\omega\mu_{\star}\varepsilon+1) + (\omega^2(1-\mu_{\star}^2\varepsilon)-i\omega\mu_{\star}-D_1f_{\pm})=0
  \end{equation}
  for some $\omega\in\R$, and for $+$ or $-$, belongs to $\sigma_{\mathrm{ess}}(\PL)$. We introduce a new unknown $\tilde{\lambda}\in\C$ 
  via $\lambda=\tilde{\lambda}+i\omega\mu_{\star}$ and solve the transformed equation
  \begin{equation*}
    \tilde{\lambda}^2 + \frac{1}{\varepsilon}\tilde{\lambda} + \frac{1}{\varepsilon}(\omega^2-D_1 f_{\pm}) = 0.
  \end{equation*}
  obtained from \eqref{equ:2.20}. Thus, the quadratic eigenvalue problem \eqref{equ:2.20} has the solutions
  \begin{align*}
    \lambda = -\frac{1}{2\varepsilon} + i\omega\mu_{\star} \pm \frac{1}{2\varepsilon}\sqrt{1-4\varepsilon(\omega^2-D_1 f_{\pm})},\,\omega\in\R.
  \end{align*}
  These solutions lie on the line $\mathrm{Re}=-\frac{1}{2\varepsilon}$ and on two ellipses if $-4D_1 f_{\pm}\varepsilon<1$ (cf. Figure \ref{fig:2.3_alt}(a)).


  \begin{figure}[H]
  \centering
  \subfigure[]{\includegraphics[height=3.9cm] {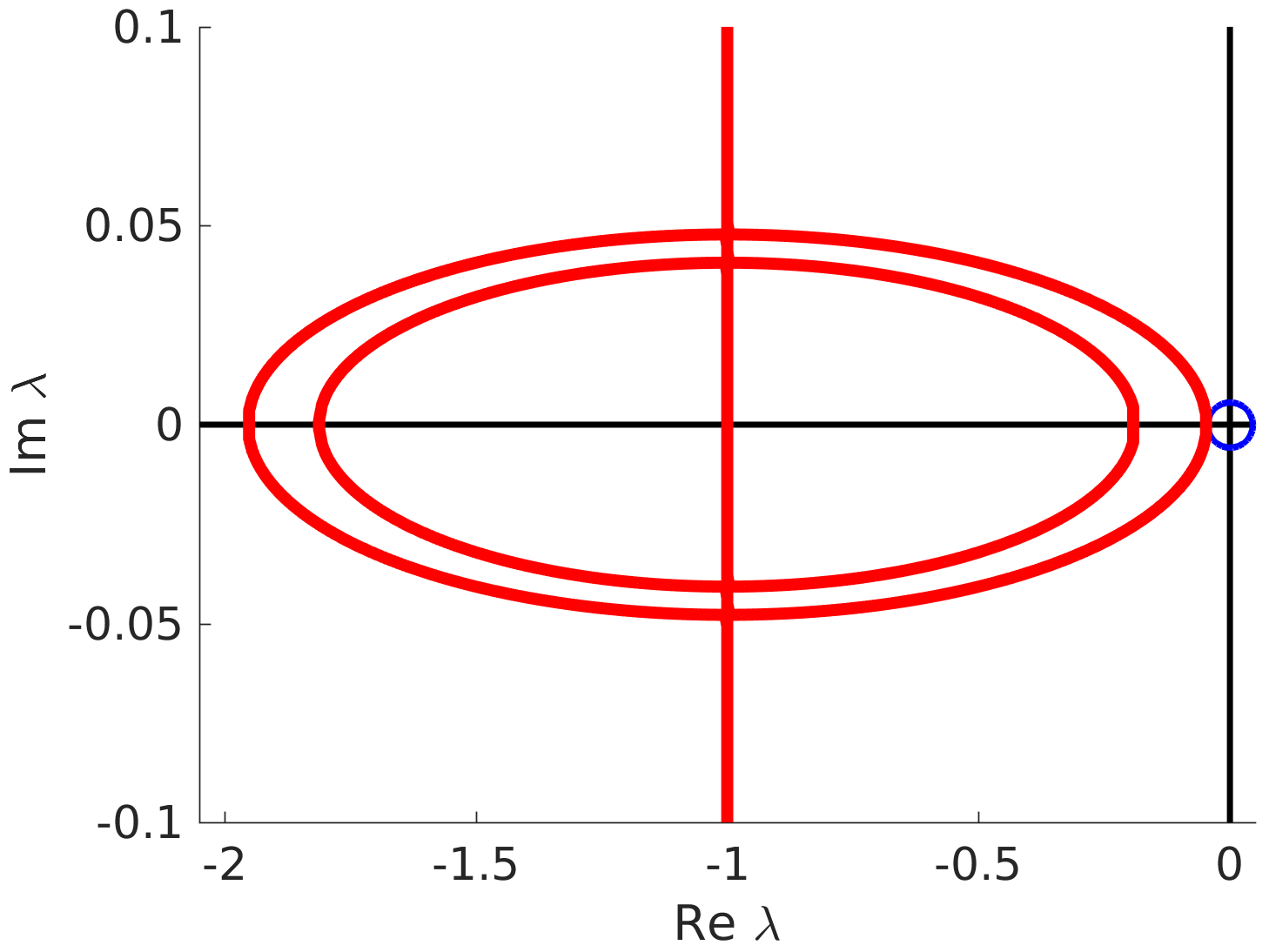}\label{fig:2.3a_alt}}
  \subfigure[]{\includegraphics[height=3.9cm] {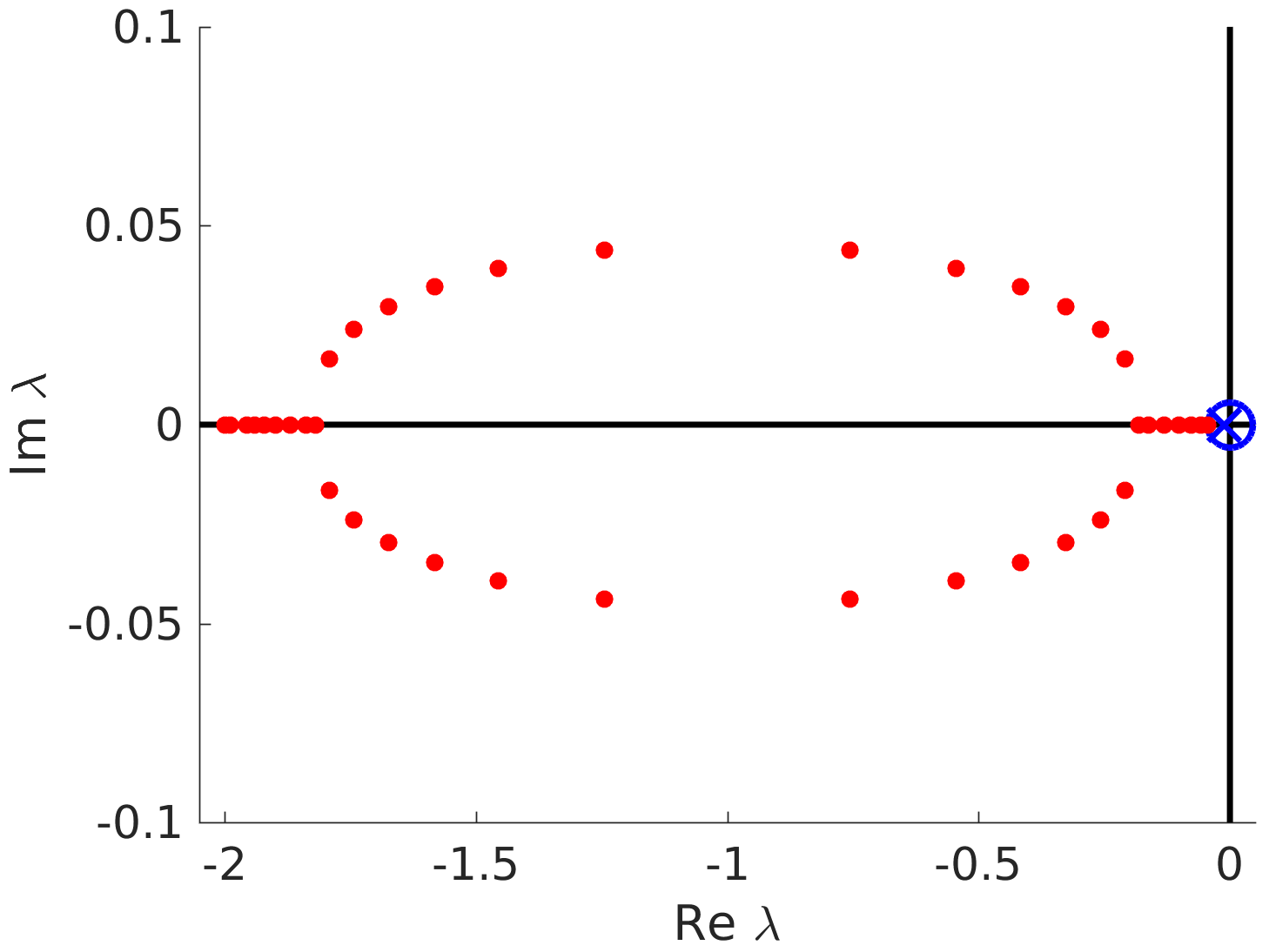} \label{fig:2.3b_alt}}
  \subfigure[]{\includegraphics[height=3.9cm] {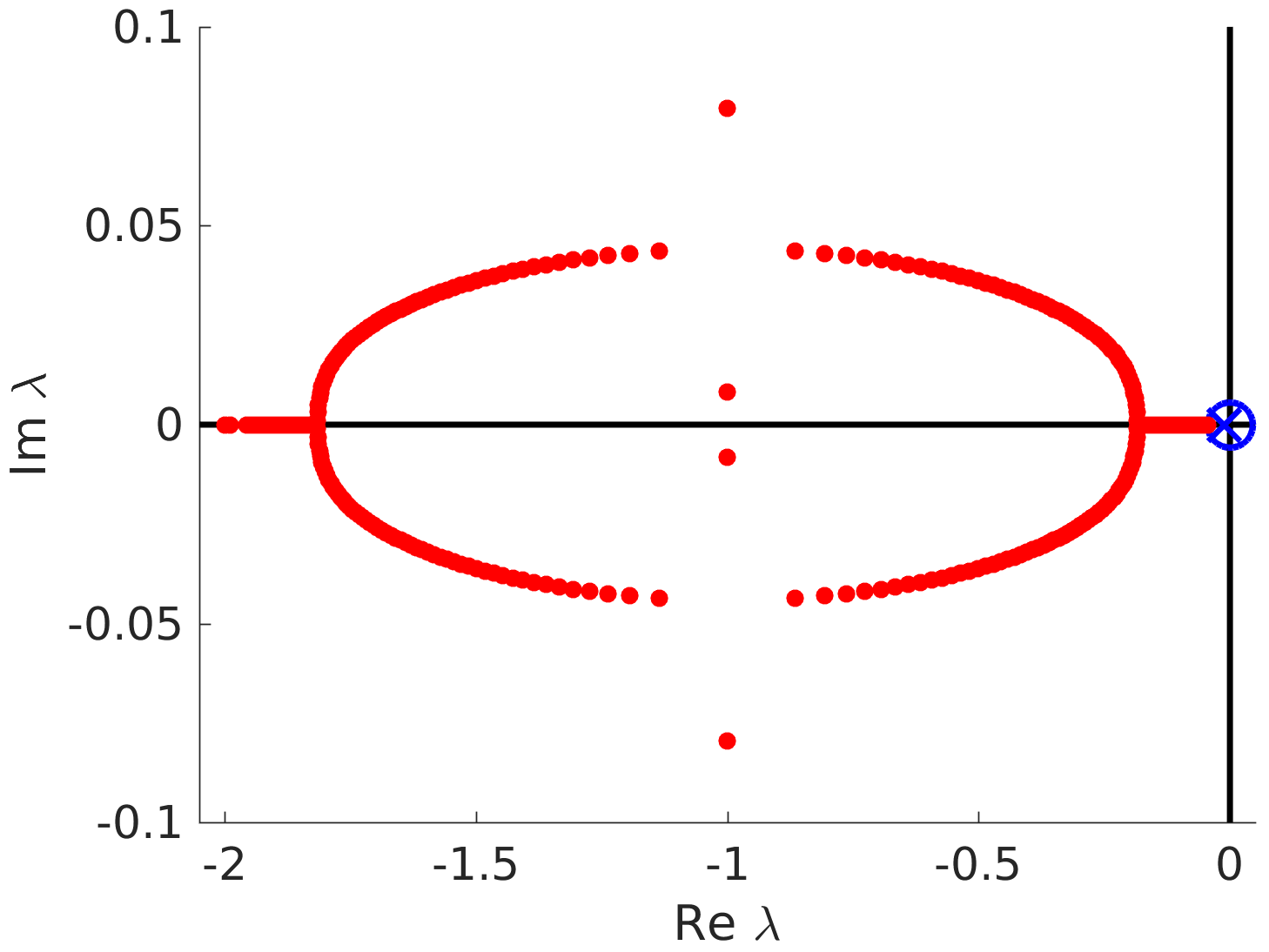}\label{fig:2.3c_alt}}
  \caption{Spectrum of the quintic Nagumo wave equation for parameters \eqref{equ:2.23} (a)
  and the numerical spectrum on the spatial domain $[-R,R]$ for $R=50$ (b) and $R=400$ (c) both for spatial stepsize $\triangle x=0.1$.}
  \label{fig:2.3_alt}
  \end{figure}

  \begin{figure}[H]
  \centering
  \subfigure[]{\includegraphics[height=3.9cm] {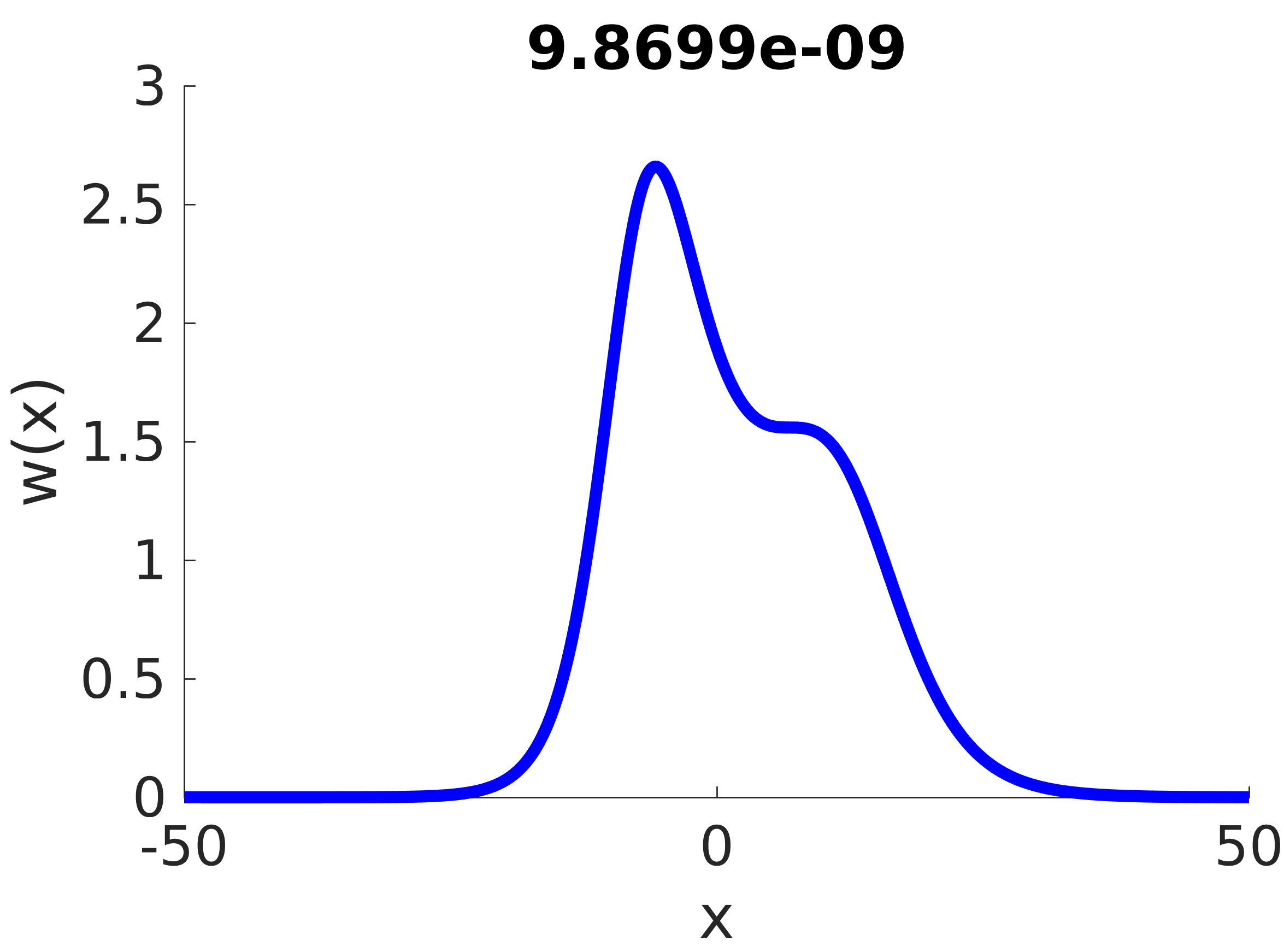}\label{fig:2.4a}}
  \subfigure[]{\includegraphics[height=3.9cm] {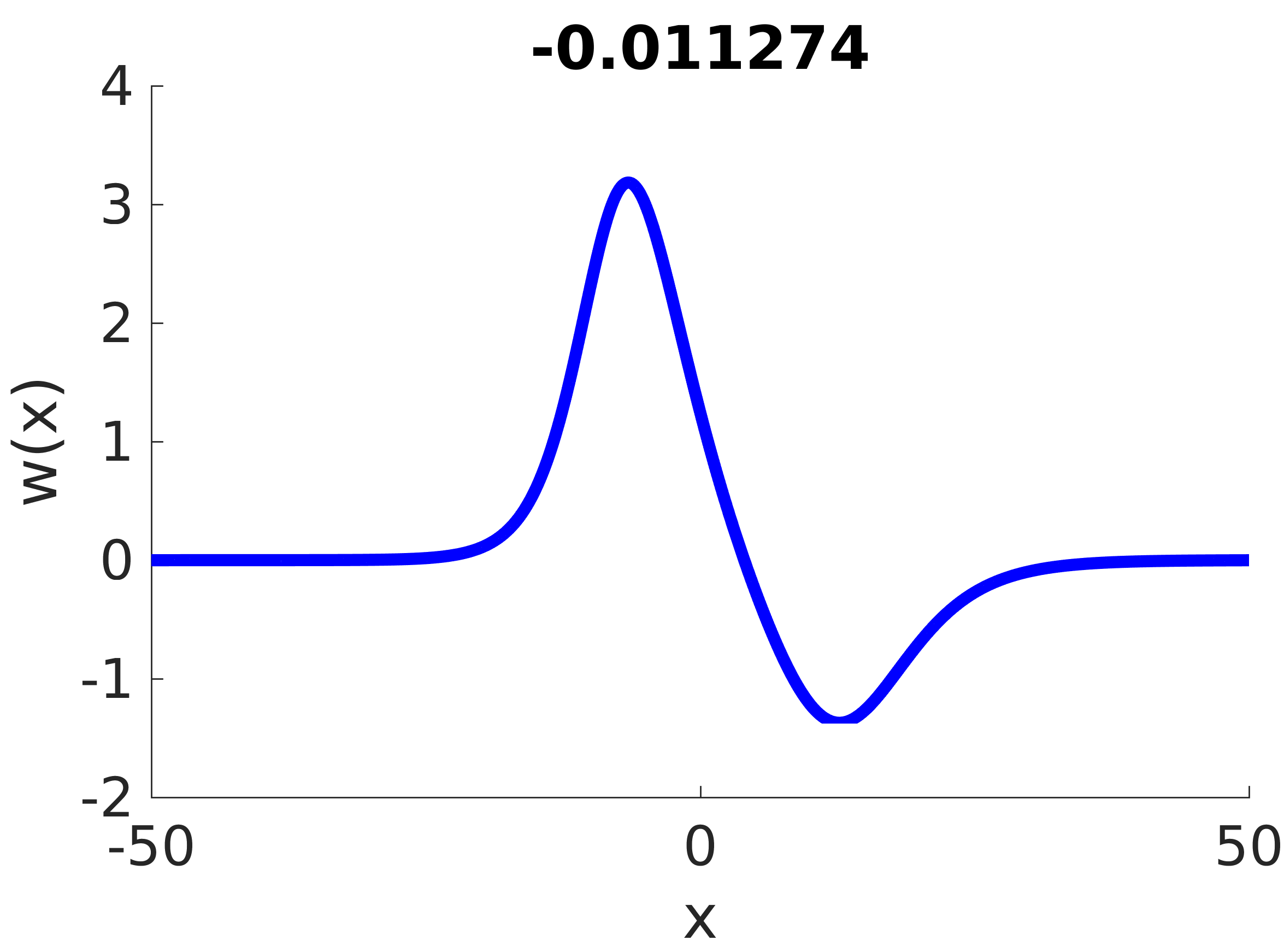} \label{fig:2.4b}}
  \subfigure[]{\includegraphics[height=3.9cm] {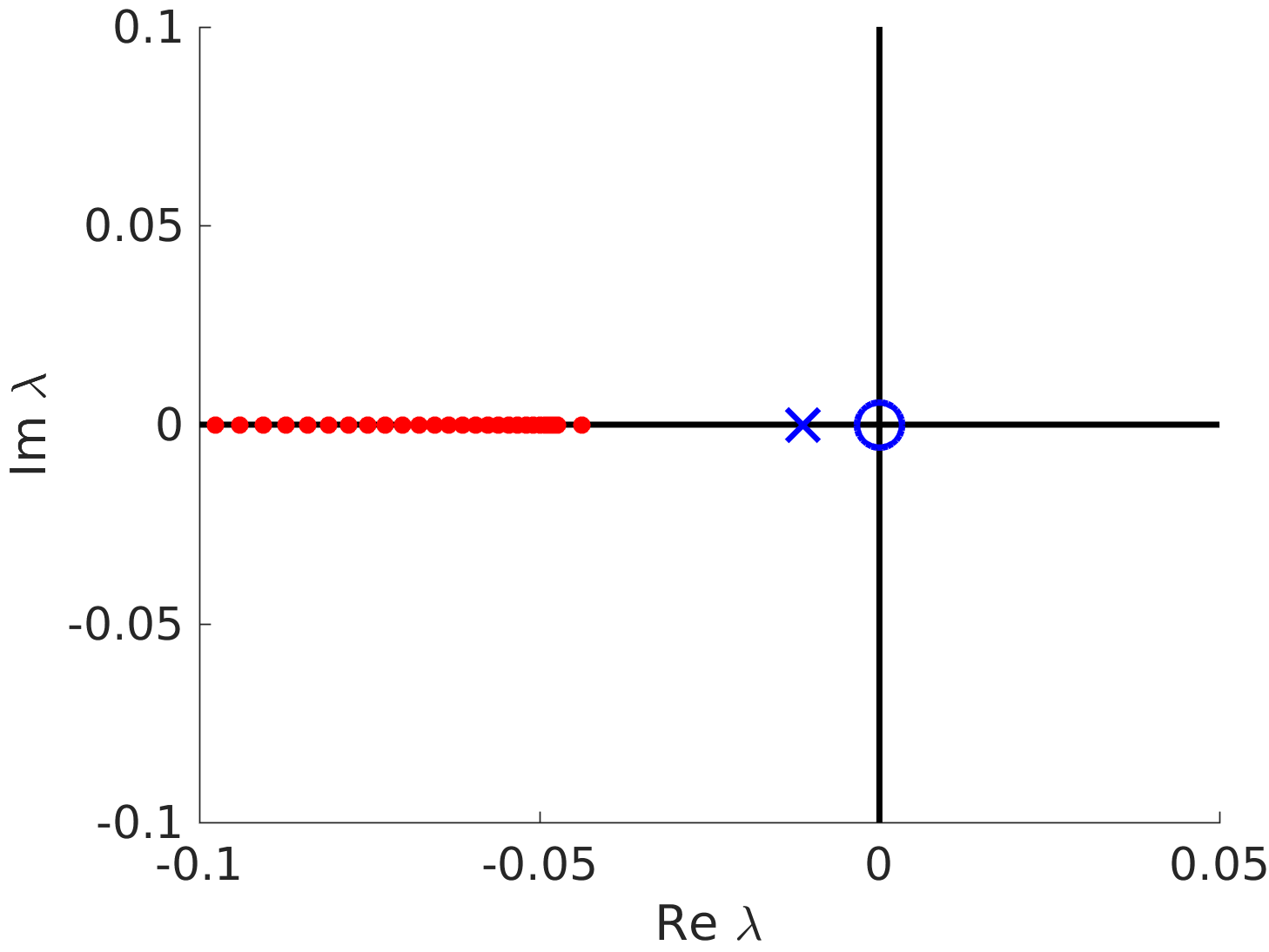} \label{fig:2.4c}}
  \caption{Eigenfunctions of the quintic Nagumo wave equation for parameters \eqref{equ:2.23} belonging to the isolated eigenvalues 
  $\lambda_1\approx 0$ (a), $\lambda_2\approx -0.011274$ (b), and a zoom into the spectrum from Fig.\ref{fig:2.3_alt}(c) in (c).}
  \label{fig:2.4}
  \end{figure}

  Figure \ref{fig:2.3_alt}(a) shows the part of the spectrum of the quintic Nagumo wave which is guaranteed by Proposition \ref{prop:2.1} and \ref{prop:2.2}. 
  It is subdivided into the symmetry set $\sigma_{\mathrm{sym}}(\PL)$ (blue circle), which is determined by Proposition \ref{prop:2.1} and belongs to 
  the point spectrum $\sigma_{\mathrm{pt}}(\PL)$, and the dispersion set $\sigma_{\mathrm{disp}}(\PL)$ (red lines), which is determined by Proposition \ref{prop:2.2} 
  and belongs to the essential spectrum $\sigma_{\mathrm{ess}}(\PL)$. In general, there may be further essential spectrum in
 $\sigma_{\mathrm{ess}}(\PL)\setminus  \sigma_{\mathrm{disp}}(\PL)$ and further isolated eigenvalues in  
  $\sigma_{\mathrm{pt}}(\PL) \setminus \sigma_{\mathrm{sym}}(\PL) $. In fact,
 for the quintic 
  Nagumo wave equation we find an extra eigenvalue with negative real part, cf. Figure \ref{fig:2.4}(c). The numerical spectrum 
  of the quintic Nagumo wave equation on the spatial domain $[-R,R]$ equipped with periodic boundary conditions is shown in Figure \ref{fig:2.3_alt}(b) for 
  $R=50$ and in Figure \ref{fig:2.3_alt}(c) for $R=400$. Each of them consists of the approximations of the point spectrum subdivided into the symmetry set 
  (blue circle) and an additional isolated eigenvalue (blue plus sign), and of the essential spectrum (red dots). 
  The missing line inside the ellipse in Figure \ref{fig:2.3_alt}(b) gradually appears numerically when enlarging the spatial domain, 
  see Figure \ref{fig:2.3_alt}(c). The second ellipse only develops on even larger domains.
\end{example}

%
%
\sect{Rotating waves in several space dimensions}
\label{sec:3}

\subsection{Freezing rotating waves.}
\label{subsec:3.1}
Consider the Cauchy problem associated with \eqref{equ:1.6}
\begin{subequations} 
  \label{equ:3.1}
  \begin{align}
    & Mu_{tt}+Bu_t = A\triangle u + f(u),         && \,x\in\R^d,\,t>0, \label{equ:3.1a} \\
    & u(\cdot,0) = u_0,\quad u_t(\cdot,0) = v_0,  && \,x\in\R^d,\,t=0, \label{equ:3.1b}
  \end{align}
\end{subequations}
for some initial data $u_0,v_0:\R^d\rightarrow\R^m$, where $u_0$ denotes the \begriff{initial displacement} and 
$v_0$ the \begriff{initial velocity}. The damped wave equation \eqref{equ:3.1} 
has a more special nonlinearity than in the one-dimensional case,
see \eqref{equ:1.6}. This will simplify some of the computations below. 

In the following, let $\SE(d)=\SO(d)\ltimes\R^d$ denote the \begriff{special Euclidean group} and $\SO(d)$ the special orthogonal group. 
Let us introduce new unknowns $(Q(t),\tau(t))\in\SE(d)$ and $v(\xi,t)\in\R^m$ via the \begriff{rotating wave ansatz}
\begin{equation}
  \begin{aligned}
  \label{equ:3.2}
    u(x,t) & = v(\xi,t),\quad\xi:=Q(t)^{\top}(x-\tau(t)),\,x\in\R^d,\,t\geqslant 0.
  \end{aligned}
\end{equation}
Inserting \eqref{equ:3.2} into \eqref{equ:3.1a} and suppressing arguments
of $u$ and $v$ leads to
\begin{align}
  \triangle_x u =& \triangle_{\xi} v,\quad f(u)=f(v), \quad
  u_t = v_{\xi}\left(Q_t^{\top}(x-\tau)-Q^{\top} \tau_t\right) + v_t,            \label{equ:3.3}\\
  u_{tt} =& v_{\xi\xi}\left(Q_t^{\top}(x-\tau)-Q^{\top} \tau_t\right)^2 + v_{\xi}\left(Q_{tt}^{\top}(x-\tau)-2Q_t^{\top} \tau_t - Q^{\top} \tau_{tt}\right) \nonumber \\
+& 2v_{\xi t}\left(Q_t^{\top}(x-\tau)-Q^{\top} \tau_t\right) + v_{tt}.
  \nonumber 
\end{align}
Hence equation \eqref{equ:3.1a} turns into
\begin{equation}
  \begin{aligned}
    \label{equ:3.5}
    \begin{split}
    &Mv_{tt} + Bv_t = A\triangle v - Mv_{\xi\xi}\left(Q_t^{\top} Q\xi - Q^{\top}\tau_t\right)^2 - 2Mv_{\xi t}\left(Q_t^{\top} Q\xi - Q^{\top}\tau_t\right) \\
    &\quad\quad\quad\quad\quad\quad- Mv_{\xi}\left(Q_{tt}^{\top} Q\xi - 2Q_t^{\top}\tau_t - Q^{\top}\tau_{tt}\right) - Bv_{\xi}\left(Q_t^{\top} Q\xi-Q^{\top}\tau_t\right) + f(v).
    \end{split}
  \end{aligned}
\end{equation}
It is convenient to introduce time-dependent functions $S_1(t),S_2(t)\in\R^{d,d}$, $\mu_1(t),\mu_2(t)\in\R^d$ via 
\begin{equation*}
  \begin{aligned}
    S_1 := Q^{\top} Q_t,\quad S_2 := S_{1,t},\quad \mu_1 := Q^{\top}\tau_t,\quad \mu_2 := \mu_{1,t}.
  \end{aligned}
\end{equation*}
Obviously, $S_1$ and $S_2$ satisfy $S_1^{\top}=-S_1$ and $S_2^{\top}=-S_2$, which follows from $Q^{\top}Q=I_d$ by differentiation.
Moreover, we obtain
\begin{align*}
  &Q_t^{\top}Q = -S_1,\quad Q^{\top}\tau_t = \mu_1,\quad Q_t^{\top}\tau_t+Q^{\top}\tau_{tt} = \mu_2, \\ 
  &Q_{tt}^{\top}Q = -S_{1,t} - S_1^{\top} S_1 = -S_2+S_1^2,\quad -Q_t^{\top}\tau_t = -Q_t^{\top} Q Q^{\top} \tau_t = S_1\mu_1,
\end{align*}
which transforms \eqref{equ:3.5} into the system
\begin{subequations} 
  \label{equ:3.6}
  \begin{align}
    &Mv_{tt} + Bv_t = A\triangle v - Mv_{\xi\xi}\left(S_1\xi+\mu_1\right)^2
 +2Mv_{\xi t}\left(S_1\xi+\mu_1\right) 
\label{equ:3.6a}\\
    &\quad\quad\quad\quad\quad\quad + Mv_{\xi}\left((S_2-S_1^2)\xi - S_1\mu_1 + \mu_2\right) + Bv_{\xi}\left(S_1\xi+\mu_1\right) + f(v), \nonumber\\
    &\begin{pmatrix}S_1\\\mu_1\end{pmatrix}_t = \begin{pmatrix}S_2\\\mu_2\end{pmatrix},\label{equ:3.6b}\\
    &\begin{pmatrix}Q\\\tau\end{pmatrix}_t = \begin{pmatrix}QS_1\\Q\mu_1\end{pmatrix}.\label{equ:3.6c}
  \end{align}
\end{subequations}
The quantity $(Q(t),\tau(t))$ describes the \begriff{position} by its spatial \begriff{shift} $\tau(t)$ and the \begriff{rotation} $Q(t)$. Moreover, $S_1(t)$ denotes the 
\begriff{rotational velocities}, $\mu_1(t)$ the \begriff{translational velocities}, $S_2(t)$ the \begriff{angular acceleration} and $\mu_2(t)$ the \begriff{translational acceleration} 
of the rotating wave $v$ at time $t$. Note that in contrast to the traveling waves the leading part $A\triangle-M\partial_{\xi}^2(\cdot)(S_1\xi+\mu_1)^2$ 
not only depends on the velocities $S_1$ and $\mu_1$, but also on the spatial variable $\xi$, which means that the leading part has unbounded (linearly growing) 
coefficients. We next specify initial data for the system \eqref{equ:3.6} as follows,
\begin{equation}
  \begin{aligned}
  \label{equ:3.7}
    \begin{split}
    &\quad v(\cdot,0) = u_0,\quad v_t(\cdot,0)=v_0+u_{0,\xi}(S_1^0\xi+\mu_1^0),\\
    &S_1(0)=S_1^0,\quad \mu_1(0)=\mu_1^0,\quad Q(0)=I_d,\quad \tau(0)=0.
    \end{split}
  \end{aligned}
\end{equation}
Note that, requiring $Q(0)=I_d$, $\tau(0)=0$, $S_1(0)=S_1^0$ and $\mu_1(0)=\mu_1^0$ for some $S_1^0\in\R^{d,d}$ with $(S_1^0)^{\top}=-S_1^0$ and $\mu_1^0\in\R^d$, 
the first equation in \eqref{equ:3.7} follows from \eqref{equ:3.2} and \eqref{equ:3.1b}, while the second condition in \eqref{equ:3.7} can be deduced from 
\eqref{equ:3.3}, \eqref{equ:3.1b}, \eqref{equ:3.6c} and the first condition in \eqref{equ:3.7}. 

The system \eqref{equ:3.6} comprises evolution equations for the unknowns $v$, $S$ and $\mu_1$.
In order to specify the remaining variables $S_2$ and $\mu_2$ we impose $\mathrm{dim}\,\SE(d)=\frac{d(d+1)}{2}$ additional scalar algebraic constraints, also 
known as \begriff{phase conditions}
\begin{align}
  \psi(v,v_t,(S_1,\mu_1),(S_2,\mu_2))=0\in\R^{\frac{d(d+1)}{2}},\quad t\geqslant 0. \label{equ:3.8} 
\end{align}
Two possible choices of such a phase condition are
\begin{align}
&\psi_{\mathrm{fix}}(v) := \begin{pmatrix}\langle v-\hat{v},D_l\hat{v}\rangle_{L^2}\\ 
                                                             \langle v-\hat{v},D^{(i,j)}\hat{v}\rangle_{L^2}\end{pmatrix} = 0,\;t\geqslant 0, \label{equ:3.9} \\
&\psi_{\mathrm{orth}}(v,v_t) := \begin{pmatrix}\langle v_t,D_l v\rangle_{L^2}\\ 
                                                             \langle v_t,D^{(i,j)} v\rangle_{L^2}\end{pmatrix} = 0,\;t\geqslant 0, \label{equ:3.10}
\end{align}
for $l=1,\ldots,d$, $i=1,\ldots,d-1$ and $j=i+1,\ldots,d$ with $D_l:=\partial_{\xi_l}$ and $D^{(i,j)}:=\xi_j \partial_{\xi_i}-\xi_i 
\partial_{\xi_j}$.
Condition \eqref{equ:3.9} is obtained from the requirement that the distance
\begin{align*}
  \rho(Q,\tau) := \left\|v(\cdot,t)-\hat{v}(Q^{\top}(\cdot -\tau))\right\|^2_{L^2}
\end{align*}
attains a local minimum at $(Q,\tau)=(I_d,0)$. Since $D_l, D^{(i,j)}$ are the generators of the Euclidean group action, condition \eqref{equ:3.10}
requires the time derivative of $v$ to be orthogonal to the group
orbit of $v$ at any time instance.

Combining the differential equations \eqref{equ:3.6}, the initial data \eqref{equ:3.7} and the phase condition \eqref{equ:3.8}, we obtain the following 
\begriff{partial differential algebraic evolution equation} (PDAE)
\begin{subequations} 
  \label{equ:3.11}
  \begin{align}
    &Mv_{tt} + Bv_t = A\triangle v - Mv_{\xi\xi}\left(S_1\xi+\mu_1\right)^2 
 + 2Mv_{\xi t}\left(S_1\xi+\mu_1\right)
\label{equ:3.11a}\\
    &\quad\quad\quad\quad\quad\quad + Mv_{\xi}\left((S_2-S_1^2)\xi - S_1\mu_1 + \mu_2\right)+ Bv_{\xi}\left(S_1\xi+\mu_1\right) + f(v), &&\,\xi\in\R^d,\,t>0,\nonumber\\
    &v(\cdot,0) = u_0,\quad v_t(\cdot,0) = v_0+u_{0,\xi}(S_{1}^0\xi+\mu_1^0), &&\,\xi\in\R^d,\,t=0, \label{equ:3.11b}\\
    &0 = \psi(v,v_t,(S_1,\mu_1),(S_2,\mu_2)),       &&\,t\geqslant 0, \label{equ:3.11c}\\ 
    &\begin{pmatrix}S_1\\\mu_1\end{pmatrix}_t = \begin{pmatrix}S_2\\\mu_2\end{pmatrix},\quad \begin{pmatrix}S_1(0)\\\mu_1(0)\end{pmatrix} = \begin{pmatrix}S_1^0\\\mu_1^0\end{pmatrix}, &&\,t\geqslant 0, \label{equ:3.11d}\\
      &\begin{pmatrix}Q\\\tau\end{pmatrix}_t = \begin{pmatrix}QS_1\\Q\mu_1\end{pmatrix},\quad \begin{pmatrix}Q(0)\\\tau(0)\end{pmatrix} = \begin{pmatrix}I_d\\0\end{pmatrix}, &&\,t\geqslant 0. 
      \label{equ:3.11e}
  \end{align}
\end{subequations}
The system \eqref{equ:3.11} depends on the choice of phase condition and must be solved for $(v,S_1,\mu_1,S_2,\mu_2,Q,\tau)$ for  given $(u_0,v_0,S_1^0,\mu_1^0)$. 
It consists of a PDE for $v$ in \eqref{equ:3.11a}--\eqref{equ:3.11b}, 
two systems of ODEs for $(S_1,\mu_1)$ in \eqref{equ:3.11d} and for $(Q,\tau)$ in \eqref{equ:3.11e} and $\frac{d(d+1)}{2}$ algebraic constraints for $(S_2,\mu_2)$ in \eqref{equ:3.11c}. 
The ODE \eqref{equ:3.11e} for $(Q,\tau)$ is the \begriff{reconstruction equation} (see \cite{RowleyKevrekidisMarsden2003}), it decouples from the other equations in \eqref{equ:3.11} 
and can be solved in a postprocessing step. Note that in the frozen equation for first order evolution equations, the ODE for $(S_1,\mu_1)$ does not appear, see \cite[(10.26)]{Otten2014}. 
The additional ODE is a new component of the PDAE and is caused by the second order time derivative.

Finally, note that $(v,S_1,\mu_1,S_2,\mu_2)=(v_{\star},S_{\star},\mu_{\star},0,0)$ satisfies
\begin{align*}
  &0 = A\triangle v - Mv_{\star,\xi\xi}\left(S_{\star}\xi+\mu_{\star}\right)^2 
        - Mv_{\star,\xi}S_{\star}\left(S_{\star}\xi + \mu_{\star}\right)+
 Bv_{\star,\xi}\left(S_{\star}\xi+\mu_{\star}\right) + f(v_{\star}),\,\xi\in\R^d,\\
  &0 = \begin{pmatrix}S_2\\\mu_2\end{pmatrix}.
 \end{align*}
 If, in addition, it has been arranged that $v_{\star},S_{\star},\mu_{\star}$ 
satisfy the phase condition $\psi(v_{\star},0,S_{\star},\mu_{\star},0,0)=0$ then  $(v_{\star},S_{\star},\mu_{\star},0,0)$ is a stationary solution of the system \eqref{equ:3.11a},\eqref{equ:3.11c},\eqref{equ:3.11d}. For a stable rotating wave we expect that solutions $(v,S_1,\mu_1,S_2,\mu_2)$ of \eqref{equ:3.11a}--\eqref{equ:3.11d} satisfy
\begin{align*}
  v(t)\rightarrow v_{\star},\quad (S_1(t),\mu_1(t))\rightarrow(S_{\star},\mu_{\star}),\quad (S_2(t),\mu_2(t))\rightarrow (0,0),\quad\text{as}\quad t\to\infty,
\end{align*}
provided the initial data are close to their limiting values.

\begin{example}[Cubic-quintic complex Ginzburg-Landau wave equation]\label{exa:3}
  Consider the cubic-quintic complex Ginzburg-Landau wave equation
  \begin{equation}
    \label{equ:3.40}
    \varepsilon u_{tt} + \rho u_t = \alpha \triangle u + u(\delta+\beta|u|^2+\gamma|u|^4),\; x\in \R^d,\, t\geqslant 0
  \end{equation}
  with $u=u(x,t)\in\C$, $d\in\{2,3\}$, $\varepsilon,\rho,\alpha,\beta,\gamma,\delta\in\C$ and $\Re\alpha>0$. 
  For the parameter values
  \begin{align}
    \label{equ:3.41}
    \varepsilon=10^{-4},\quad\rho=1,\quad\alpha=\frac{3}{5},\quad\gamma=-1-\frac{1}{10}i,\quad\beta=\frac{5}{2}+i,\quad\delta=-0.73.
  \end{align}
  equation \eqref{equ:3.40} admits a spinning soliton solution.
  
  \begin{figure}[ht]
    \centering
    \subfigure[]{\includegraphics[height=4.0cm] {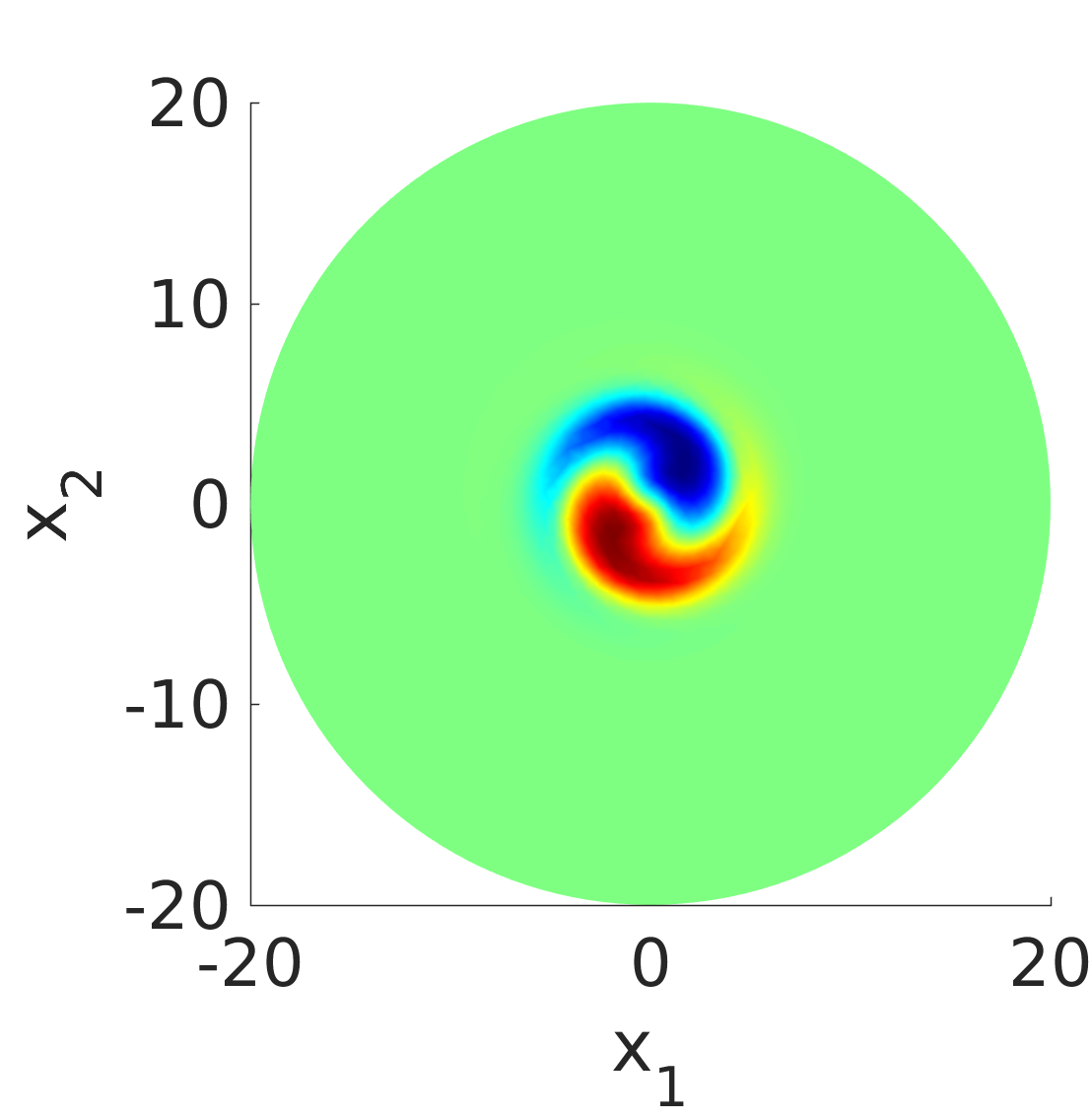}\label{fig:3.1a}}
    \subfigure[]{\includegraphics[height=4.0cm] {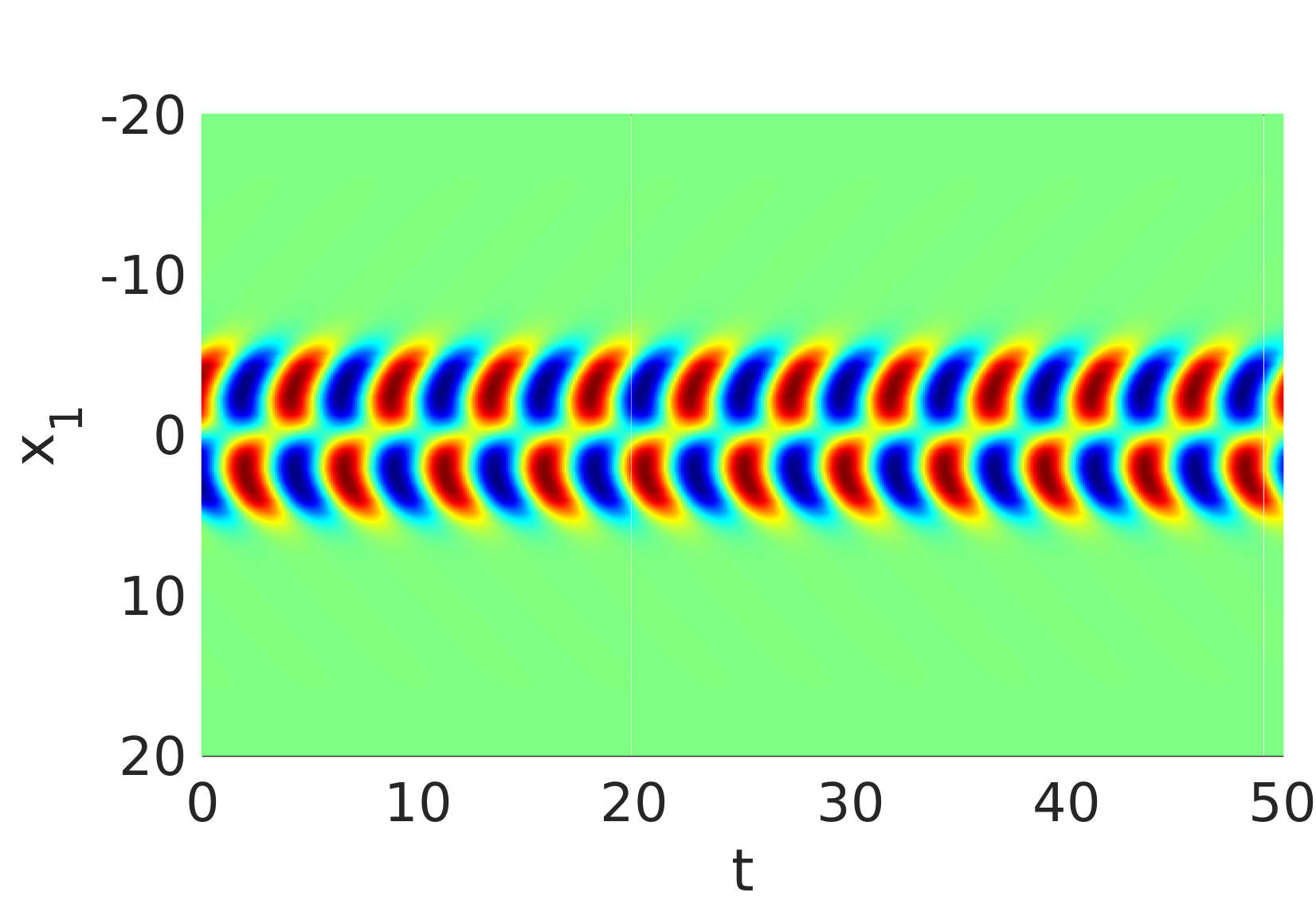} \label{fig:3.1b}}
    \caption{Solution of cubic-quintic complex Ginzburg-Landau wave equation \eqref{equ:3.40}: Spinning soliton $u(x,t)$ at time $t=50$ (a) and its time evolution along $x_2=0$ (b) 
    for parameters from \eqref{equ:3.41}.}
    \label{fig:3.1}
  \end{figure}

  Figure \ref{fig:3.1} shows a numerical simulation of the solution $u$ of \eqref{equ:3.40} on the ball $B_R(0)$ of radius $R=20$, 
  with homogeneous Neumann boundary conditions and with parameter
  values from \eqref{equ:3.41}. The initial data $u_0$ and $v_0$ are generated in the following way. First we use the freezing method
  to compute a rotating wave in the parabolic case (as in  \cite{Otten2014}) for parameter values $\varepsilon=0$, $\rho=1$ and
  \begin{align*}
   \alpha=\frac{1}{2}+\frac{1}{2}i,\quad\gamma=-1-\frac{1}{10}i,
\quad\beta=\frac{5}{2}+i,\quad\delta=-\frac{1}{2}.
  \end{align*}
  Then the parameter set $(\varepsilon,\alpha,\delta)$ is gradually
  changed until the values \eqref{equ:3.41} are attained.
   For the space discretization we use continuous piecewise 
  linear finite elements with spatial stepsize $\triangle x=0.8$. For the time discretization we use the BDF method of order $2$ 
  with absolute tolerance $\mathrm{atol}=10^{-4}$, relative tolerance $\mathrm{rtol}=10^{-3}$, temporal stepsize $\triangle t=0.1$
  and final time $T=50$. Computations are performed with the help of the software COMSOL 5.2.
   
  Let us now consider the frozen cubic-quintic complex Ginzburg-Landau wave equation resulting from \eqref{equ:3.11}
  \begin{subequations} 
  \label{equ:3.42}
  \begin{align}
    &\varepsilon v_{tt} + \rho v_t = \alpha\triangle v - \varepsilon v_{\xi\xi}\left(S_1\xi+\mu_1\right)^2  + 2\varepsilon v_{\xi t}\left(S_1\xi+\mu_1\right) \label{equ:3.42a}\\
    &\quad\quad\quad\quad\quad\quad + \varepsilon v_{\xi}\left((S_2-S_1^2)\xi - S_1\mu_1 + \mu_2\right)+ \rho v_{\xi}\left(S_1\xi+\mu_1\right) + f(v), &&\,\xi\in\R^d,\,t>0,\nonumber\\
    &v(\cdot,0) = u_0,\quad v_t(\cdot,0) = v_0+u_{0,\xi}(S_{1}^0\xi+\mu_1^0), &&\,\xi\in\R^d,\,t=0, \label{equ:3.42b}\\
    &0 = \psi_{\mathrm{fix}}(v) := \begin{pmatrix}\langle v-\hat{v},D_l\hat{v}\rangle_{L^2}\\ 
                                                   \langle v-\hat{v},D^{(i,j)}\hat{v}\rangle_{L^2}\end{pmatrix}, &&\,t\geqslant 0, \label{equ:3.42c}\\ 
    &\begin{pmatrix}S_1\\\mu_1\end{pmatrix}_t = \begin{pmatrix}S_2\\\mu_2\end{pmatrix},\quad 
     \begin{pmatrix}S_1(0)\\\mu_1(0)\end{pmatrix} = \begin{pmatrix}S_1^0\\\mu_1^0\end{pmatrix}, &&\,t\geqslant 0, \label{equ:3.42d}\\
    &\begin{pmatrix}Q\\\tau\end{pmatrix}_t = \begin{pmatrix}QS_1\\Q\mu_1\end{pmatrix},\quad 
     \begin{pmatrix}Q(0)\\\tau(0)\end{pmatrix} = \begin{pmatrix}I_d\\0\end{pmatrix}, &&\,t\geqslant 0. \label{equ:3.42e}
  \end{align}
  \end{subequations}

  \begin{figure}[ht]
    \centering
    \subfigure[]{\includegraphics[height=4.0cm] {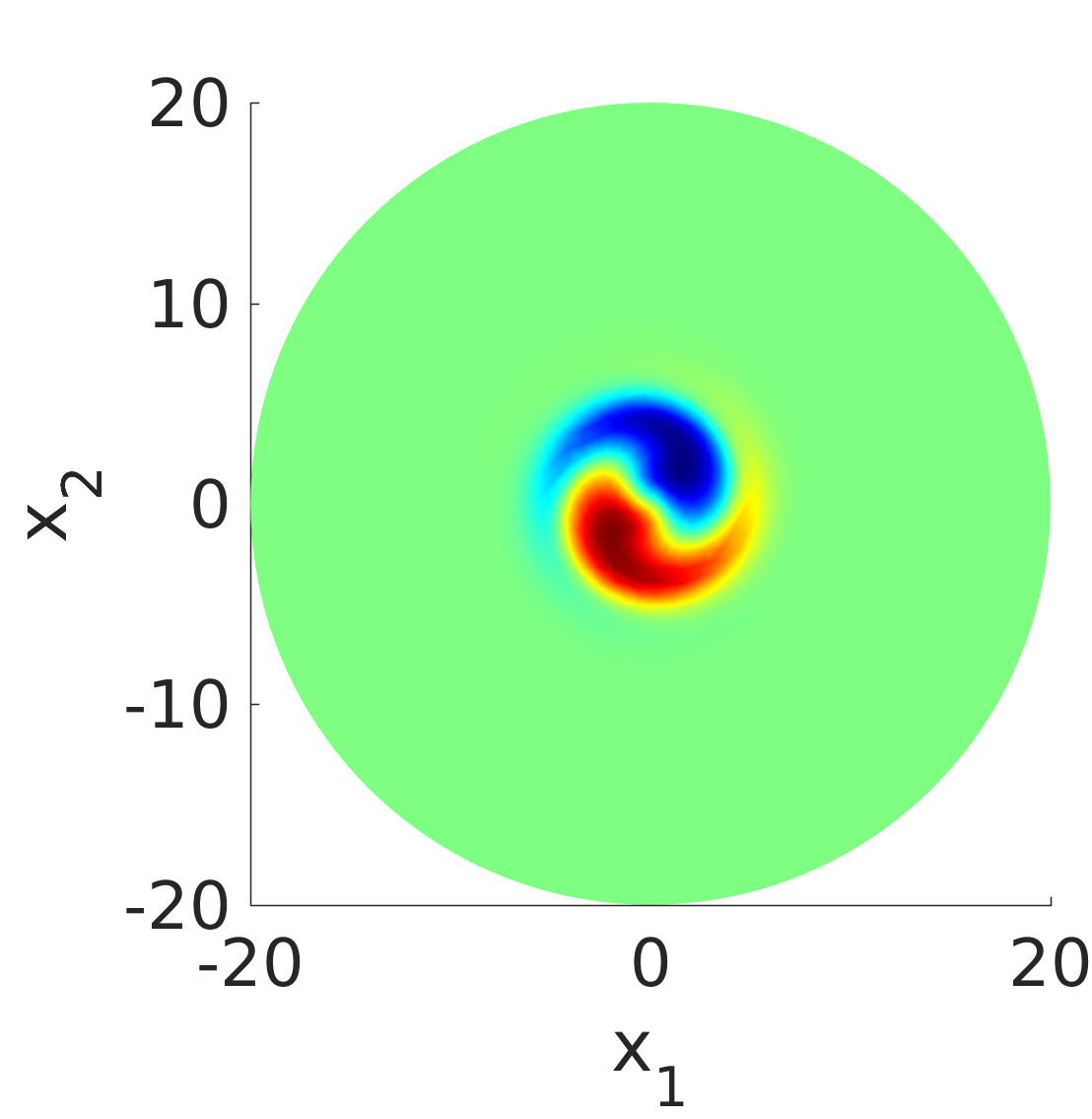}\label{fig:3.2a}}
    \subfigure[]{\includegraphics[height=4.0cm] {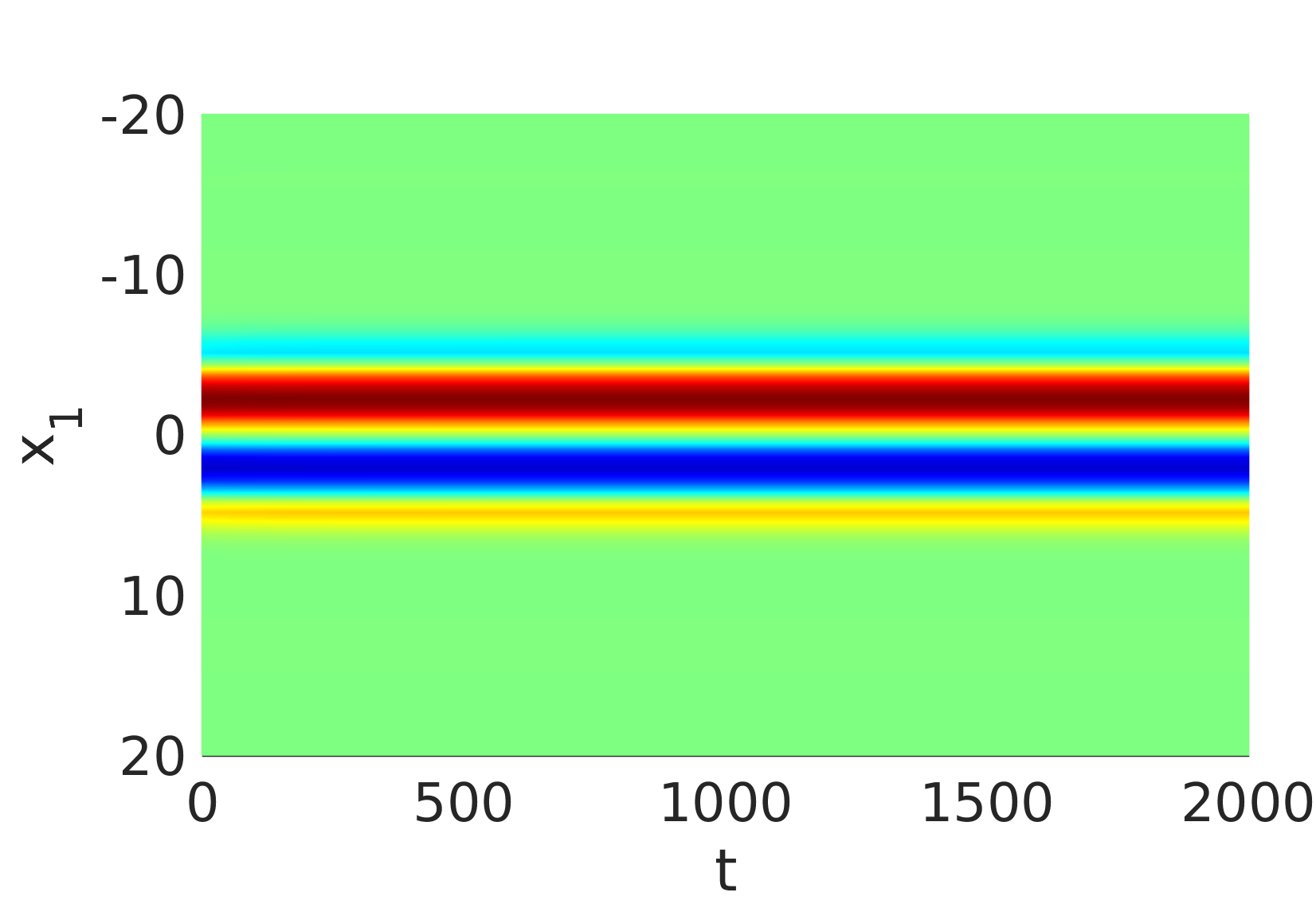} \label{fig:3.2b}}\\
    \subfigure[]{\includegraphics[height=4.0cm] {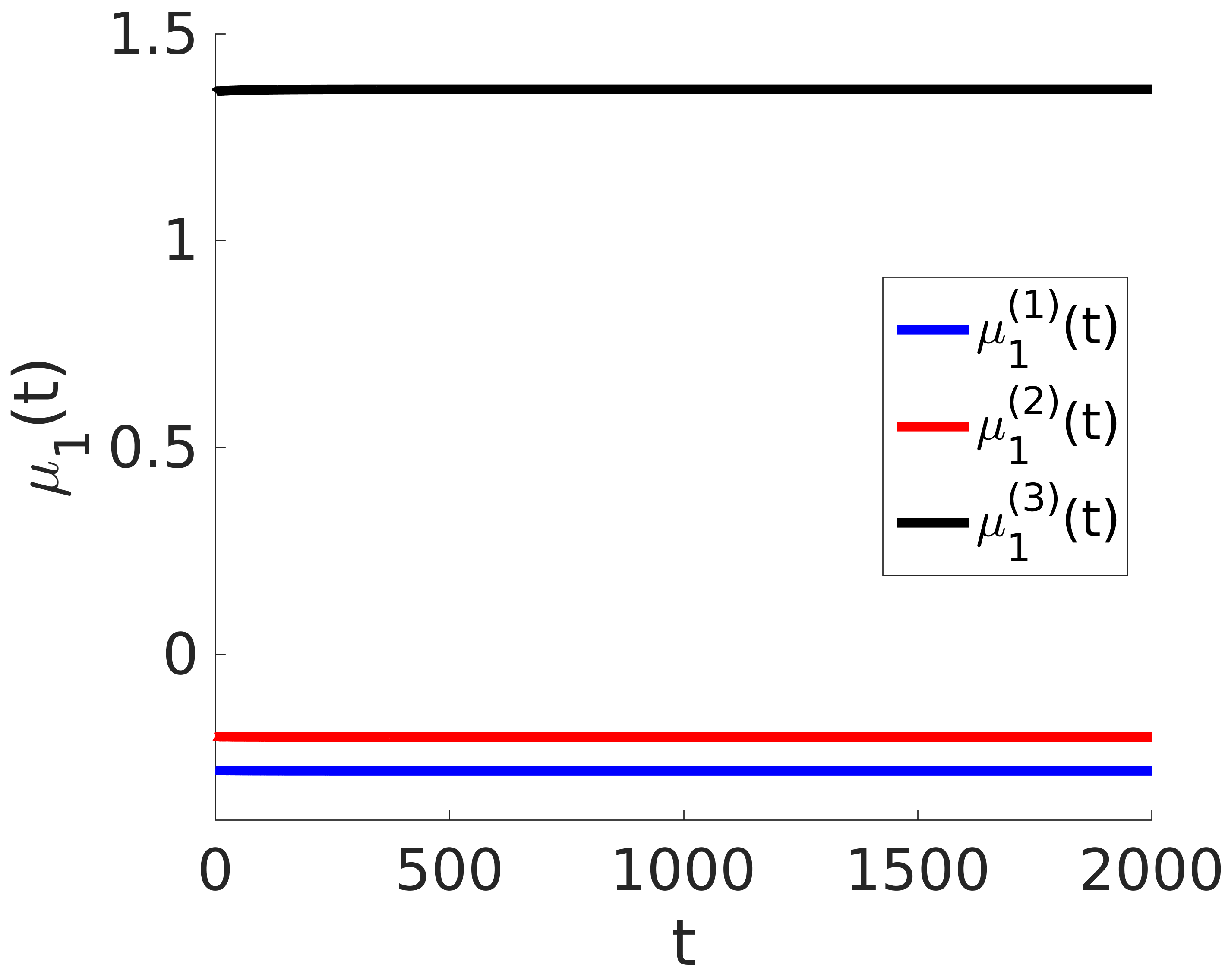}\label{fig:3.2c}}
    \subfigure[]{\includegraphics[height=4.0cm] {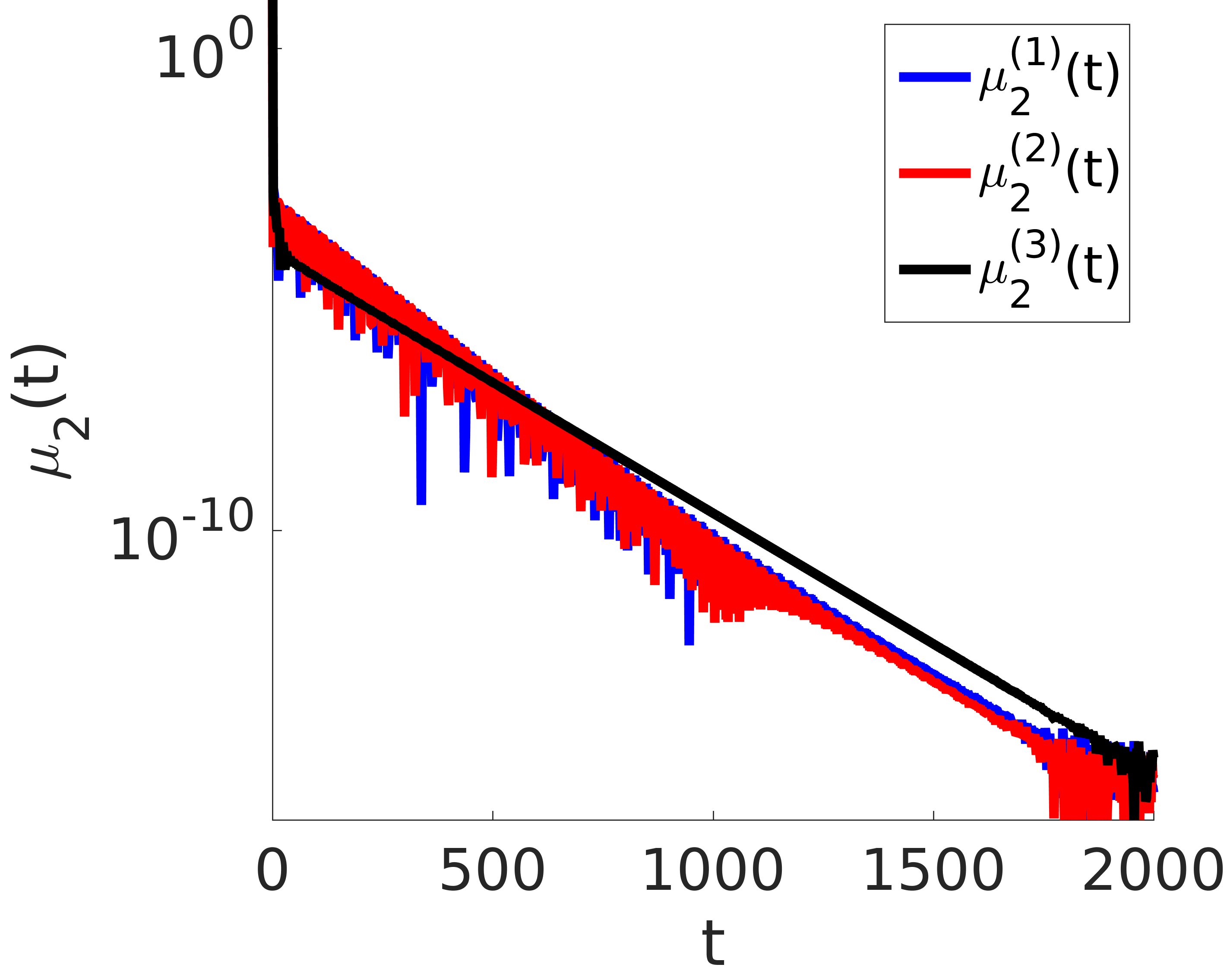} \label{fig:3.2d}}
    \caption{Solution of the frozen cubic-quintic complex Ginzburg-Landau wave equation \eqref{equ:3.42}: profile $v(x,t)$ at time $t=2000$ (a), 
    its time evolution along $x_2=0$ (b), velocities $\mu_1(t)$ (c), and accelerations $\mu_2(t)$ (d) for parameters from \eqref{equ:3.41}.}
    \label{fig:3.2}
  \end{figure}

  Figure \ref{fig:3.2} shows the solution $(v,S_1,\mu_1,S_2,\mu_2,Q,\tau)$ of \eqref{equ:3.42} on the ball $B_R(0)$ with radius $R=20$, 
  homogeneous Neumann boundary conditions, initial data $u_0$, $v_0$ as in the nonfrozen case, and reference function $\hat{v}=u_0$. 
  For the computation we used the fixed phase condition $\psi_{\mathrm{fix}}(v)$ from \eqref{equ:3.9}. The spatial discretization data are 
  taken as in the nonfrozen case. For the time discretization we used the BDF method of order $2$ with absolute tolerance $\mathrm{atol}=10^{-3}$, 
  relative tolerance $\mathrm{rtol}=10^{-2}$, maximal temporal stepsize $\triangle t=0.5$, initial step $10^{-4}$, and final time $T=2000$. 
  Due to the choice of initial data, the profile becomes immediately stationary, the acceleration $\mu_2$ converges to zero, while
  the speed $\mu_1$ and the nontrivial entry $S_{12}$ of $S$ approach asymptotic values 
  \begin{align*}
    \mu_1^{(1)}=-0.2819,\quad\mu_1^{(2)}=-0.1999,\quad S_{12}=1.3658.
  \end{align*}
  Note that we have a clockwise rotation if $S_{12}>0$, and a counter clockwise rotation if $S_{12}<0$. Thus, 
  the spinning soliton rotates clockwise. The center of rotation $x_{\star}$ and the temporal period $T^{\mathrm{2D}}$ for one rotation are given by, see \cite[Exa.10.8]{Otten2014},
  \begin{align*}
    x_{\star} = \frac{1}{S_{12}}\begin{pmatrix}\mu_1^{(2)}\\-\mu_1^{(1)}\end{pmatrix}=\begin{pmatrix}-0.1464\\0.2064\end{pmatrix},\quad\quad
    T^{\mathrm{2D}}=\frac{2\pi}{|S_{12}|}=4.6004.
  \end{align*}
\end{example}

\subsection{Spectra of rotating waves.}
\label{subsec:3.2}

Consider the linearized equation
\begin{equation}
    \label{equ:3.13}
      Mv_{tt} + Bv_t - A\triangle v + Mv_{\xi\xi}(S_{\star}\xi)^2 - 2Mv_{\xi t}S_{\star}\xi
      + Mv_{\xi}S_{\star}^2\xi
               - Bv_{\xi}S_{\star}\xi  
                - Df(v_{\star})v = 0.
\end{equation}
Equation \eqref{equ:3.13} is obtained from
the co-rotating frame equation \eqref{equ:1.8} when linearizing at the profile $v_{\star}$. Moreover, we assume $\mu_{\star}=0$, that is the wave 
rotates about the origin. Shifting the center of rotation does not influence the stability properties, see the discussion in \cite{BeynOtten2016}. 
Looking for solutions of the form $v(\xi,t)=e^{\lambda t}w(\xi)$ to \eqref{equ:3.13} yields the quadratic eigenvalue problem
\begin{equation}
  \begin{aligned}
    \label{equ:3.14}
    \PL(\lambda)w := \left(\lambda^2 P_2 + \lambda P_1 + P_0\right)w = 0,\,\xi\in\R^d
  \end{aligned}
\end{equation}
with differential operators $P_j$ defined by
\begin{equation}
  \begin{aligned}
    \label{equ:3.15}
    \begin{split}
    P_2 =& M,\quad 
    P_1 = B-2M\left(\partial_{\xi}\,\cdot\right)S_{\star}\xi = B-
    2M\sum_{j=1}^{d}(S_{\star}\xi)_j\partial_{\xi_j}, \\
    P_0 =& -A\triangle\,\cdot + M\left(\partial_{\xi}^2\,\cdot\right)(S_{\star}\xi)^2
            + M\left(\partial_{\xi}\,\cdot\right)S_{\star}^2\xi 
-B\left(\partial_{\xi}\,\cdot\right)S_{\star}\xi- Df(v_{\star})\,\cdot \\
        =& -A\sum_{j=1}^{d}\partial_{\xi_j}^2 + 
M\sum_{j=1}^{d}\sum_{\nu=1}^{d}(S_{\star}\xi)_j(S_{\star}\xi)_\nu\partial_{\xi_j}\partial_{\xi_\nu} + 
M\sum_{j=1}^{d}(S_{\star}^2\xi)_j\partial_{\xi_j}
            - B\sum_{j=1}^{d}(S_{\star}\xi)_j \partial_{\xi_j} - Df(v_{\star}).
    \end{split}
  \end{aligned}
\end{equation}
As in the one-dimensional case we cannot solve equation \eqref{equ:3.14} in general. Rather, our aim is to determine the 
symmetry set $\sigma_{\mathrm{sym}}(\PL)$ as a subset of the point spectrum $\sigma_{\mathrm{pt}}(\PL)$, and the dispersion 
set $\sigma_{\mathrm{disp}}(\PL)$ as a subset of the essential spectrum $\sigma_{\mathrm{ess}}(\PL)$. The point spectrum 
is affected by the underlying group symmetries while the essential spectrum depends on the far-field behavior of the wave.

In the following we present the recipe for computing the subsets
$\sigma_{\mathrm{sym}}(\PL)\subseteq \sigma_{\mathrm{pt}}(\PL)$ and
$\sigma_{\mathrm{disp}}(\PL)\subseteq \sigma_{\mathrm{ess}}(\PL)$.

\subsubsection{Point Spectrum and symmetry set.}
\label{subsubsec:3.2.2}
Let us look for eigenfunctions $w$ of \eqref{equ:3.14} of the form
\begin{align}
  \label{equ:3.26}
  w(\xi) = v_{\star,\xi}(\xi)(E\xi+b)\quad\quad\text{for some $E\in\C^{d,d}$, $b\in\C^d$, $E^{\top}=-E$, $v_{\star}\in C^3(\R^d,\R^m)$.}
\end{align}
This ansatz is motivated by the fact that functions of this type
span the image of the derivative of the
group action $(Q,\tau)\rightarrow v_{\star}(Q^{\top}(\cdot - \tau))$ at
the unit element $(Q,\tau)=(I_d,0)\in \SE(d)$ (compare \eqref{equ:3.2}).
We plug \eqref{equ:3.26} into \eqref{equ:3.14} and use the equalities 
\begin{align}
  &Mw = Mv_{\star,\xi}(E\xi+b),\quad\quad
   Bw = Bv_{\star,\xi}(E\xi+b), \nonumber\\
  &2M(\partial_{\xi} w)S_{\star}\xi = 2M v_{\star,\xi\xi}(E\xi+b)S_{\star}\xi + 2Mv_{\star,\xi}ES_{\star}\xi
  \label{equ:3.27} \\
  &A\triangle w = (\partial_{\xi}(A\triangle v_{\star}))(E\xi+b)
  \label{equ:3.28} \\
  &M(\partial_{\xi}^2 w)(S_{\star}\xi)^2 = (\partial_{\xi}(M v_{\star,\xi\xi}(S_{\star}\xi)^2))(E\xi+b) + 2Mv_{\star,\xi\xi}([E,S_{\star}]\xi-S_{\star}b)S_{\star}\xi,
  \label{equ:3.29} \\
  &M(\partial_{\xi} w)S_{\star}^2\xi = (\partial_{\xi}(M v_{\star,\xi}S_{\star}^2\xi))(E\xi+b) + Mv_{\star,\xi}([E,S_{\star}^2]\xi-S_{\star}^2 b),
  \label{equ:3.30} \\
  &B(\partial_{\xi} w)S_{\star}\xi = (\partial_{\xi}(Bv_{\star,\xi}S_{\star}\xi))(E\xi+b) + Bv_{\star,\xi}([E,S_{\star}]\xi-S_{\star}b),
  \label{equ:3.31} \\
  &Df(v_{\star})w = (\partial_{\xi}(f(v_{\star})))(E\xi+b)
  \label{equ:3.32}
\end{align}
where $[E,S_{\star}]:=ES_{\star}-S_{\star}E$ is the Lie bracket.
This leads to the following equation:
\begin{align}
  \label{equ:3.33}
  0 =& \lambda^2 M v_{\star,\xi}(E\xi+b) + \lambda\Big(Bv_{\star,\xi}(E\xi+b) - 2Mv_{\star,\xi\xi}(E\xi+b)S_{\star}\xi - 2Mv_{\star,\xi}ES_{\star}\xi\Big) \nonumber\\
     & +\Big(2Mv_{\star,\xi\xi}([E,S_{\star}]\xi-S_{\star}b)S_{\star}\xi + Mv_{\star,\xi}([E,S_{\star}^2]\xi-S_{\star}^2 b) -Bv_{\star,\xi}([E,S_{\star}]\xi-S_{\star}b) \\
     &\;\;\quad\;-\partial_{\xi}\big(A\triangle v_{\star} - Mv_{\star,\xi\xi}(S_{\star}\xi)^2 - Mv_{\star,\xi}S_{\star}^2\xi + Bv_{\star,\xi}S_{\star}\xi + f(v_{\star})\big)(E\xi+b)\Big). \nonumber
\end{align}
Now we use the rotating wave equation \eqref{equ:1.9} in 
\eqref{equ:3.33} and obtain by rearranging the remaining terms
\begin{align}
  \label{equ:3.35}
  0 =& Mv_{\star,\xi}\Big(\lambda^2(E\xi+b) - 2\lambda ES_{\star}\xi + [E,S_{\star}^2]\xi - S_{\star}^2 b\Big)
       + Bv_{\star,\xi}\Big(\lambda(E\xi+b)-[E,S_{\star}]\xi+S_{\star}b\Big) \nonumber \\
     & - 2Mv_{\star,\xi\xi}\Big(\lambda(E\xi+b) - [E,S_{\star}]\xi + S_{\star}b\Big)S_{\star}\xi \\
    =& Mv_{\star,\xi}\Big( (\lambda^2 E - 2\lambda ES_{\star} + [E,S_{\star}^2])\xi + \lambda^2 b - S_{\star}^2 b\Big) 
       + Bv_{\star,\xi}\Big((\lambda E-[E,S_{\star}])\xi+\lambda b+S_{\star}b\Big) \nonumber \\
     & - 2Mv_{\star,\xi\xi}\Big((\lambda E-[E,S_{\star}])\xi + \lambda b+S_{\star}b\Big)S_{\star}\xi. \nonumber
\end{align}
Comparing coefficients in \eqref{equ:3.35} yields the finite-dimensional eigenvalue problem (see \cite{BlochIserles2005},\cite{Otten2014}, \cite{BeynOtten2016b})
\begin{subequations}
  \label{equ:3.36}
  \begin{align}
  \lambda E &= [E,S_{\star}],
  \label{equ:3.36a} \\
  \lambda b &= -S_{\star}b,
  \label{equ:3.36b}  
  \end{align}
\end{subequations}
which must be solved for $(\lambda,E,b)$ and admits $\frac{d(d+1)}{2}$ solutions. In fact, having a solution $(\lambda,E,b)$ of \eqref{equ:3.36}, then 
the last two terms in \eqref{equ:3.35} obviously vanish. The first term vanishes if we write both summands as
\begin{align*}
  \lambda^2 b - S_{\star}^2 b = \lambda(\lambda b+S_{\star}b)-S_{\star}(\lambda b+S_{\star}b)
\end{align*}
and
\begin{align*}
  &\lambda^2 E - 2\lambda ES_{\star} + [E,S_{\star}^2]
  = \lambda(\lambda E-[E,S_{\star}])-(2\lambda ES_{\star}-\lambda[E,S_{\star}]-[E,S_{\star}^2]) \\
  =& \lambda(\lambda E-[E,S_{\star}])-\left((\lambda E-[E,S_{\star}])S_{\star}+S_{\star}(\lambda E-[E,S_{\star}])+[E,S_{\star}]S_{\star}+S_{\star}[E,S_{\star}]-[E,S_{\star}^2]\right),
\end{align*}  
and use the identity $[E,S_{\star}]S_{\star}+S_{\star}[E,S_{\star}]-[E,S_{\star}^2]=[E,[S_{\star},S_{\star}]]=0$ which holds by skew-symmetry of $S_{\star}$. 
Therefore, it is sufficient to solve \eqref{equ:3.36}. Furthermore, if $(\lambda,E)$ is a solution of \eqref{equ:3.36a}, then $(\lambda,E,0)$ solves \eqref{equ:3.36}, and,
similarly, if $(\lambda,b)$ is a solution 
of \eqref{equ:3.36b}, then $(\lambda,0,b)$ solves \eqref{equ:3.36}. Therefore, it is sufficient to solve \eqref{equ:3.36a} and \eqref{equ:3.36b} separately. 
For the skew-symmetric matrix $S_{\star}$  we have $S_{\star}=U\Lambda U^{\herm}$ for some unitary $U \in C^{d,d}$ and some diagonal matrix
$\Lambda=\diag(\lambda_1,\ldots,\lambda_d)$ where
$\lambda_1,\ldots,\lambda_d \in i \R$ are the eigenvalues of $S_{\star}$. 
In particular, this implies $S_{\star}^{\top}=\overline{U}\Lambda U^{\top}$.

\begin{itemize}[leftmargin=*]
\item Multiply \eqref{equ:3.36b} from the left by $U^{\herm}$ and define $\tilde{b}=U^{\herm}b$ to obtain
\begin{align}
  \label{equ:3.37}
  \lambda \tilde{b} = \lambda U^{\herm}b = -U^{\herm}S_{\star}b = -U^{\herm}U\Lambda U^{\herm}b = -\Lambda \tilde{b}.
\end{align}
Equation \eqref{equ:3.37} has solutions $(\lambda,\tilde{b})=(-\lambda_l,e_l)$, hence \eqref{equ:3.36b} has solutions $(\lambda,b)=(-\lambda_l,Ue_l)$, 
and \eqref{equ:3.36} has solutions $(\lambda,E,b)=(-\lambda_l,0,Ue_l)$ for $l=1,\ldots,d$.

\item Multiply \eqref{equ:3.36a} from the left by $U^{\herm}$, from
  the right by $\bar{U}$, define $\tilde{E}=U^{\herm}E\overline{U}$, 
and use the skew-symmetry of $S_{\star}$ and $\tilde{E}$, to obtain
\begin{align}
  \label{equ:3.38}
  \lambda\tilde{E} = \lambda U^{\herm}E\overline{U} = U^{\herm}[E,S_{\star}]\overline{U} = -U^{\herm}E\overline{U}\Lambda U^{\top}\overline{U} - U^{\herm}U\Lambda U^{\herm}E\overline{U} 
  = -\tilde{E}\Lambda - \Lambda \tilde{E} = \tilde{E}^{\top}\Lambda - \Lambda\tilde{E}.
\end{align}
Equation \eqref{equ:3.38} has solutions $(\lambda,\tilde{E})=(-(\lambda_i+\lambda_j),I_{ij}-I_{ji})$, hence \eqref{equ:3.36a} has solutions 
$(\lambda,E)=(-(\lambda_i+\lambda_j),U(I_{ij}-I_{ji})U^{\top})$, and \eqref{equ:3.36} has solutions $(\lambda,E,b)=(-(\lambda_i+\lambda_j),U(I_{ij}-I_{ji})U^{\top},0)$ 
for $i=1,\ldots,d-1$, $j=i+1,\ldots,d$, where $I_{ij}$ has entry $1$ in the $i$th row and $j$th column and $0$ otherwise.
\end{itemize}

Let us summarize the result in a proposition.

\begin{proposition}[Point spectrum of rotating waves]\label{prop:3.2}
  Let $f\in C^2(\R^m,\R^m)$  and let $v_{\star}\in C^3(\R^d,\R^m)$ be a classical solution of 
  \eqref{equ:1.9} for some skew-symmetric matrix $S_{\star}\in \R^{m,m}$ with eigenvalues denoted 
  by $\lambda_1,\ldots,\lambda_d$ and unitary matrix $U\in\C^{d,d}$ diagonalizing $S_{\star}$. Then
  \begin{align*}
    w=v_{\star,\xi}(E\xi+b)
  \end{align*}
  is a classical solution of the eigenvalue problem \eqref{equ:3.14} provided that
  \begin{align*}
    (\lambda,E,b)=(-\lambda_l,0,Ue_l)\quad\quad\text{or}\quad\quad (\lambda,E,b)=(-(\lambda_i+\lambda_j),U(I_{ij}-I_{ji})U^{\top},0)
  \end{align*}
  for some $l=1,\dots,d$, $i=1,\ldots,d-1$, $j=i+1,\ldots,d$. In particular, the symmetry set
  \begin{align*}
    \sigma_{\mathrm{sym}}(\PL) = \sigma(S_{\star}) \cup \{\lambda_i+\lambda_j:1\leqslant i<j\leqslant d\}.
  \end{align*}
  belongs to the point spectrum $\sigma_{\mathrm{pt}}(\PL)$ of $\PL$.
\end{proposition}

Altogether, Proposition \ref{prop:3.2} yields $\frac{d(d+1)}{2}$ solutions of the quadratic eigenvalue problem \eqref{equ:3.14}. 
It is a remarkable feature that the eigenvalues and the eigenfunctions coincide with those for first order evolution equations, 
see \cite{BeynOtten2016b}, \cite{Otten2014}. Moreover, we suggest that Proposition \ref{prop:3.2} also applies to rotating waves 
that are not localized, e.g. spiral waves and scroll waves. This has been confirmed in numerical experiments. 

\begin{figure}[ht]
  \centering
  \subfigure[$d=2$\newline{}$\mathrm{dim}\,\SE(2)=3$]{\includegraphics[page=1, height=6cm]{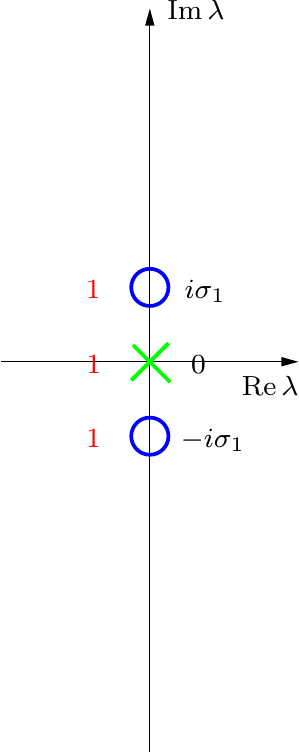} \label{fig:Pointspectrum_d2}}
  \subfigure[$d=3$\newline{}$\mathrm{dim}\,\SE(3)=6$]{\includegraphics[page=2, height=6cm]{Images.pdf} \label{fig:Pointspectrum_d3}}
  \subfigure[$d=4$\newline{}$\mathrm{dim}\,\SE(4)=10$]{\includegraphics[page=3, height=6cm]{Images.pdf} \label{fig:Pointspectrum_d4}}
  \subfigure[$d=5$\newline{}$\mathrm{dim}\,\SE(5)=15$]{\includegraphics[page=4, height=6cm]{Images.pdf}\label{fig:Pointspectrum_d5}}
  \caption{Point spectrum of the linearization $\PL$ on the imaginary axis $i\R$ for space dimension $d=2,3,4,5$ given by Proposition \ref{prop:3.2}.}
  \label{fig:3.5}
\end{figure}

Figure \ref{fig:3.5} shows the eigenvalues $\lambda\in\sigma_{\mathrm{sym}}(\PL)$ from Proposition \ref{prop:3.2} and their corresponding 
multiplicities for different space dimensions $d=2,3,4,5$. The eigenvalues $\lambda\in\sigma(S_{\star})$ are indicated by blue circles, 
the eigenvalues $\lambda\in\left\{\lambda_i+\lambda_j\mid \lambda_i,\lambda_j\in\sigma(S_{\star}),\,1\leqslant i<j\leqslant d\right\}$ 
by green crosses. The imaginary values to the right of the symbols denote eigenvalues and the numbers to the left their corresponding 
multiplicities. As expected, there are $\frac{d(d+1)}{2}$ eigenvalues on the imaginary axis in case of space dimension $d$. 

\subsubsection{Essential spectrum and dispersion set.}\label{subsubsec:3.2.1}

\begin{enumerate}[label=\bf{\arabic*.},leftmargin=*]
\item \textbf{Quasi-diagonal real form.}
Let us transform the skew-symmetric matrix $S_{\star}$ into quasi-diagonal real form. For this purpose, let $\pm i\sigma_1,\ldots,\pm i\sigma_k$ be the nonzero eigenvalues 
of $S_{\star}$ so that $0$ is a semisimple eigenvalue of multiplicity $d-2k$. There is an orthogonal matrix $P\in\R^{d,d}$ such that
\begin{equation*} 
  S_{\star}=P\Lambda P^{\top},\quad
  \Lambda=\mathrm{diag}\left(\Lambda_1,\ldots,\Lambda_k,\mathbf{0}\right),\quad
  \Lambda_j=\begin{pmatrix}0 &\sigma_j\\-\sigma_j &0\end{pmatrix},\quad
  \mathbf{0}\in\R^{d-2k,d-2k}.
\end{equation*}
The transformation $\tilde{w}(y)=w(Py),\tilde{v}_{\star}(y)=v_{\star}(Py)$ transfers \eqref{equ:3.14} with operators $P_j$ from \eqref{equ:3.15} into
\begin{equation} \label{equ:3.14a}
(\lambda^2\tilde{P}_2+ \lambda \tilde{P}_1 + \tilde{P}_0)\tilde{w} =0.
\end{equation}
With the abbreviations
\begin{equation} \label{equ:3.17a}
 D_j=\partial_{y_j}, \quad D^{(i,j)}=y_jD_i - y_i D_j,\quad
K=\sum_{l=1}^k \sigma_l D^{(2l-1,2l)}
\end{equation}
the operators $\tilde{P}_j$ are given by
\begin{equation}
\label{equ:3.18}
  \begin{aligned}
          \tilde{P}_2 = & M,\quad
      \tilde{P}_1 =   B-2M\sum_{j=1}^{d}(\Lambda y)_jD_j
                  =   B-2MK,\\
                  \tilde{P}_{0} = & -A\triangle
                  +M\sum_{j=1}^{d}\sum_{\nu=1}^{d}(\Lambda y)_j(\Lambda y)_{\nu}
                  D_jD_{\nu}\,+ M\sum_{j=1}^{d}(\Lambda^2 y)_j D_j\,
                    - B\sum_{j=1}^{d}(\Lambda y)_j D_j\,- Df(\tilde{v}_{\star})\\
                   = & -A\triangle +
                     MK^2 - BK - Df(\tilde{v}_{\star}).
      \end{aligned}
\end{equation}
\item \textbf{The far-field operator.} Assume that $v_{\star}$ has an
asymptotic state  $v_{\infty}\in\R^m$, i.e.  $f(v_{\infty})=0$ and $v_{\star}(\xi)\to v_{\infty}\in\R^m$ as $|\xi|\to\infty$. 
In the limit $|y| \rightarrow \infty$ the eigenvalue problem 
\eqref{equ:3.14a} turns into the far-field problem
\begin{equation}
  \label{equ:3.17}
  \left(\lambda^2 \tilde{P}_2 + \lambda \tilde{P}_1 + \tilde{P}_{\infty}  
\right)\tilde{w} = 0,\,y\in\R^d, \quad \tilde{P}_{\infty}=-A\triangle +
  MK^2 - BK - Df(v_{\infty}).
\end{equation}

\item \textbf{Transformation into several planar polar coordinates.}
Since we have $k$ angular derivatives in $k$ different planes it is advisable to transform into several planar polar coordinates via
\begin{equation*}
  \begin{pmatrix}y_{2l-1}\\y_{2l}\end{pmatrix} = T(r_l,\phi_l):=\begin{pmatrix}r_l\cos\phi_l\\r_l\sin\phi_l\end{pmatrix},\;\phi_l\in[-\pi,\pi),\;r_l\in(0,\infty),\;l=1,\ldots,k.
\end{equation*}
All further coordinates, i.e. $y_{2k+1},\ldots,y_d$, remain fixed. The transformation
$\hat{w}(\psi):=\tilde{w}(T_2(\psi))$ with  $T_2(\psi) = (T(r_1,\phi_1),\ldots,T(r_k,\phi_k),y_{2k+1},\ldots,y_d)$ 
for $\psi=(r_1,\phi_1,\ldots,r_k,\phi_k,y_{2k+1},\ldots,y_d)$ in
the domain $\Omega=((0,\infty)\times[-\pi,\pi))^k\times\R^{d-2k}$ transfers \eqref{equ:3.17} into
\begin{equation}
  \label{equ:3.19}
  \left(\lambda^2 \hat{P}_2 + \lambda \hat{P}_1 + \hat{P}_{\infty}  
\right)\hat{w} = 0,\,\psi\in\Omega
\end{equation}
with
\begin{equation*}
  \begin{aligned}
    \begin{split}
      &\hat{P}_2 =  M,\quad
      \hat{P}_1 =   B+2M\sum_{l=1}^{k}\sigma_l\partial_{\phi_l}, \\
      &\hat{P}_{\infty} =  - A\bigg[\sum_{l=1}^{k}\left(\partial_{r_l}^2+\frac{1}{r_l}\partial_{\phi_l}+\frac{1}{r_l^2}\partial_{\phi_l}^2\right)+
        \sum_{l=2k+1}^{d}\partial_{y_l}^2\bigg] +
M\sum_{l,n=1}^{k}\sigma_l\sigma_n\partial_{\phi_l}\partial_{\phi_n} 
                   + B\sum_{l=1}^{k}\sigma_l\partial_{\phi_l} - Df(v_{\infty}).
      \end{split}
  \end{aligned}
\end{equation*}

\item \textbf{Simplified far-field operator:}
The far-field operator \eqref{equ:3.19} can be further simplified by letting 
$r_l\to\infty$ for any $1\leqslant l\leqslant k$ which turns \eqref{equ:3.19} into
\begin{equation}
  \label{equ:3.21}
  \left(\lambda^2 \hat{P}_2 + \lambda \hat{P}_1 +
  P_{\infty}^{\mathrm{sim}}\right)\hat{w} = 0,\,\psi\in\Omega
\end{equation}
with 
\begin{equation}
      \label{equ:3.22}
        P_{\infty}^{\mathrm{sim}} =  -A\left[\sum_{l=1}^{k}\partial_{r_l}^2+
      \sum_{l=2k+1}^{d}\partial_{y_l}^2\right]
+M\sum_{l,n=1}^{k}\sigma_l\sigma_n\partial_{\phi_l}
\partial_{\phi_n} 
                  + B\sum_{l=1}^{k}\sigma_l\partial_{\phi_l} - Df(v_{\infty}).
    \end{equation}

\item \textbf{Angular Fourier transform:}
Finally, we solve for eigenvalues and eigenfunctions of \eqref{equ:3.22}
by separation of variables and an angular resp. radial Fourier ansatz 
with $\omega\in\R^k$, $\rho,y\in\R^{d-2k}$, $n\in\Z^k$, $z\in\C^m$, $|z|=1$, $r\in(0,\infty)^k$, $\phi\in(-\pi,\pi]^k$:
\begin{align*}
  \hat{w}(\psi) = \exp\left(i\sum_{l=1}^{k}\omega_l r_l\right)\exp\left(i\sum_{l=1}^{k}n_l\phi_l\right)\exp\left(i\sum_{l=2k+1}^{d}\rho_l y_l\right)z 
                = \exp\left(i\langle\omega,r\rangle + i\langle n,\phi\rangle + i\langle\rho,y\rangle\right)z.
\end{align*}
Inserting this in \eqref{equ:3.21} leads to the $m$-dimensional quadratic eigenvalue problem
\begin{equation}
  \label{equ:3.23}
  \left(\lambda^2 A_2 + \lambda A_1(n) +
  A_{\infty}(\omega,n,\rho)\right)z = 0
\end{equation}
with matrices $A_2\in\R^{m,m}$ and $A_1(n),
A_{\infty}(\omega,n,\rho)\in\C^{m,m}$ given by
\begin{equation}
  \begin{aligned}
    \label{equ:3.24}
    \begin{split}
      A_2 = & M,\quad
      A_1(n) = B+2i\langle\sigma,n\rangle M,\\
      A_{\infty}(\omega,n,\rho) = &    \left(|\omega|^2+|\rho|^2\right)A - \langle\sigma,n\rangle^2 M + i\langle\sigma,n\rangle B - Df(v_{\infty}).
      \end{split}
  \end{aligned}
\end{equation}
The Fourier ansatz is a well-known tool for investigating essential spectra, see e.g. \cite{FiedlerScheel2003}.

\item \textbf{Dispersion relation and dispersion set:} As in Section
  \ref{subsubsec:2.2.2} we consider the dispersion set consisting of all values $\lambda \in \C$
  satisfying the dispersion relation
  \begin{equation}
  \label{equ:3.25}
  \det\left(\lambda^2 A_2 + \lambda A_1(n) +
   A_{\infty}(\omega,n,\rho)\right) = 0
  \end{equation}
  for some $\omega\in\R^k$, $\rho\in\R^{d-2k}$ and $n\in\Z^k$. Of course,
  one can replace $|\omega|^2+|\rho|^2$ by any nonnegative real number.
  Solving \eqref{equ:3.25} is equivalent to finding all zeros of a parameterized polynomial of degree $2m$. 
  Note that the limiting case $M=0$ and $B=I_m$ in \eqref{equ:3.25} leads to the dispersion relation for rotating 
  waves of first order evolution equations, see \cite{BeynLorenz2008} for $d=2$, and \cite[Sec. 7.4 and 9.4]{Otten2014}, 
  \cite{BeynOtten2016b} for general $d\geqslant 2$. 
\end{enumerate}
Using standard cut-off arguments as in \cite{BeynLorenz2008},\cite{Otten2014},\cite{BeynOtten2016b}, 
the following result can be shown for suitable function spaces (e.g. $L^2(\R^d,\R^m)$):

\begin{proposition}[Essential spectrum of rotating waves]\label{prop:3.1}
  Let $f\in C^1(\R^m,\R^m)$ with $f(v_{\infty})=0$ for some $v_{\infty}\in\R^m$. Let $v_{\star}\in C^2(\R^d,\R^m)$ with skew-symmetric $S_{\star}\in \R^{m,m}$ be a classical solution of 
  \eqref{equ:1.9} satisfying $v_{\star}(\xi)\to v_{\infty}$ as $|\xi|\to\infty$. Then, the dispersion set
  \begin{equation*}
  \sigma_{\mathrm{disp}}(\PL) = \{\lambda\in\C\mid\text{$\lambda$ satisfies \eqref{equ:3.25} for some $\omega\in\R^k$, $\rho\in\R^{d-2k}$, $n\in\Z^k$}\}
  \end{equation*}
  belongs to the essential spectrum $\sigma_{\mathrm{ess}}(\PL)$ of the operator 
polynomial $\PL$ from \eqref{equ:3.14}.
\end{proposition}

\begin{example}[Cubic-quintic Ginzburg-Landau wave equation]\label{exa:4}
  As shown in Example \ref{exa:3} the cubic-quintic Ginzburg-Landau wave equation \eqref{equ:3.40} with coefficients and parameters \eqref{equ:3.41}
  has a spinning soliton solution $u_{\star}(x,t)=v_{\star}(e^{-tS_{\star}}(x-x_{\star}))$ with rotational velocity $\mu_1^{(3)}=1.3658$.

  We next solve numerically the eigenvalue problem for the cubic-quintic Ginzburg-Landau wave equation. For this purpose we consider the real valued version 
  of \eqref{equ:3.40}
  \begin{equation}
    \label{equ:3.44}
    M U_{tt} + B U_t = A \triangle U + F(U),\; x\in \R^d,\, t\geqslant 0
  \end{equation}
  with
  \begin{align}
    \label{equ:3.45}
    \begin{split}
    &M = \begin{pmatrix}\varepsilon_1&-\varepsilon_2\\\varepsilon_2&\varepsilon_1\end{pmatrix},\quad B = \begin{pmatrix}\rho_1&-\rho_2\\\rho_2&\rho_1\end{pmatrix},\quad
    A = \begin{pmatrix}\alpha_1&-\alpha_2\\\alpha_2&\alpha_1\end{pmatrix},\quad U = \begin{pmatrix}u_1\\u_2\end{pmatrix},\\
    &F(U) = \begin{pmatrix}(U_1\delta_1-U_2\delta_2)+(U_1\beta_1-U_2\beta_2)(U_1^2+U_2^2)+(U_1\gamma_1-U_2\gamma_2)(U_1^2+U_2^2)^2\\
                          (U_1\delta_2+U_2\delta_1)+(U_1\beta_2+U_2\beta_1)(U_1^2+U_2^2)+(U_1\gamma_2+U_2\gamma_1)(U_1^2+U_2^2)^2\end{pmatrix},
    \end{split}
  \end{align}
  where $u=u_1+iu_2$, $\varepsilon=\varepsilon_1+i\varepsilon_2$, $\rho=\rho_1+i\rho_2$, $\alpha=\alpha_1+i\alpha_2$, $\beta=\beta_1+i\beta_2$, $\gamma=\gamma_1+i\gamma_2$, $\delta=\delta_1+i\delta_2$ 
  and $\varepsilon_j,\rho_j,\alpha_j,\beta_j,\gamma_j,\delta_j\in\R$. 

  Now, the eigenvalue problem for the cubic-quintic Ginzburg-Landau wave equation is, cf. \eqref{equ:3.14}, \eqref{equ:3.15},
  \begin{align}
  \label{equ:3.46}
    \left(\lambda^2 M\cdot + \lambda\left[B\cdot-2M(\partial_\xi \cdot)S\xi\right] + \left[-A\triangle\cdot + M(\partial_{\xi}^2 \cdot)(S\xi)^2+M(\partial_{\xi}\cdot)S^2\xi-B(\partial_{\xi}\cdot)S\xi-DF(v_{\star})\cdot\right]\right)w= 0.
  \end{align}
  Both approximations of the profile $v_{\star}$ and the velocity matrix $S=S_{\star}$ in \eqref{equ:3.46} are chosen from the solution of \eqref{equ:3.42} at time 
  $t=2000$ in Example \ref{exa:3}. By Proposition \ref{prop:3.2} the problem \eqref{equ:3.46} has eigenvalues $\lambda=0,\pm i\sigma$. These eigenvalues will
be isolated and hence belong to the point spectrum, if the differential
operator is Fredholm of index $0$ in suitable function spaces. For the parabolic
case ($M=0$) this has been established in \cite{BeynOtten2016b} and we expect
it to hold in the general case as well.
  Let us next discuss the dispersion set from Proposition \ref{prop:3.1}. The cubic-quintic Ginzburg-Landau nonlinearity $F:\R^2\rightarrow\R^2$ from \eqref{equ:3.45}
  satisfies
  \begin{equation}
    \label{equ:3.47}
    DF(v_{\infty})=\begin{pmatrix}\delta_1&-\delta_2\\\delta_2&\delta_1\end{pmatrix}\quad\text{for}\quad v_{\infty}=\begin{pmatrix}0\\0\end{pmatrix}.
  \end{equation}
  The matrices $A_2$, $A_1(n)$, $A_{\infty}(\omega,n)$ from \eqref{equ:3.24} of the quadratic problem \eqref{equ:3.23} are given by
  \begin{align*}
    A_2 = M,\quad
    A_1(n) = B+2i\sigma n M,\quad
    A_{\infty}(\omega,n) = \omega^2 A - \sigma^2 n^2 M + i\sigma n B - DF(v_{\infty})
  \end{align*}
  for $M,B,A$ from \eqref{equ:3.45}, $DF(v_{\infty})$ from \eqref{equ:3.47}, $\omega\in\R$, $n\in\Z$ and $\sigma=\mu_{1}^{(3)}$.
  The dispersion relation \eqref{equ:3.25} for the spinning solitons of the Ginzburg-Landau wave equation in $\R^2$ states that every $\lambda\in\C$ satisfying
  \begin{equation*}
    \det\left(\lambda^2 M + \lambda(B+2i\sigma n M) + (\omega^2 A - \sigma^2 n^2 M + i\sigma n B - DF(v_{\infty}))\right)=0
  \end{equation*}
  for some $\omega\in\R$ and $n\in\Z$, belongs to the essential spectrum $\sigma_{\mathrm{ess}}(\PL)$ of $\PL$. We may rewrite this in complex notation and find 
   the dispersion set
  \begin{align}
    \label{equ:3.50} \sigma_{\mathrm{disp}}(\PL)= \{ \lambda \in \C:
    \lambda^2 \varepsilon + \lambda(\rho+2i\sigma n \varepsilon) + (\omega^2 \alpha - \sigma^2 n^2 \varepsilon + i\sigma n \rho - \delta)=0 \; \text{for some} \;  \omega\in \R,n\in \Z\}.
  \end{align}
  The elements of the dispersion set are
  \begin{align*}
    \lambda_{1,2} = -\frac{\rho}{2\varepsilon} - i\sigma n \pm \frac{1}{2\varepsilon}\sqrt{\rho^2-4\varepsilon(\omega^2\alpha-\delta)},\quad
   n \in \Z, \omega \in \R. 
  \end{align*}
  They lie on the vertical line $\mathrm{Re}=-\frac{\rho}{2\varepsilon}$ and on infinitely many horizontal lines given for $n \in \Z$ by \\$i\sigma n + \frac{1}{2 \varepsilon}
  [ -\rho - \sqrt{\rho^2+4 \varepsilon \delta}, \rho+
    \sqrt{\rho^2+4 \varepsilon \delta}]$,
see Figure \ref{fig:3.3} (a),(b).

 \begin{figure}[H]
  \centering
  \subfigure[]{\includegraphics[height=4.2cm] {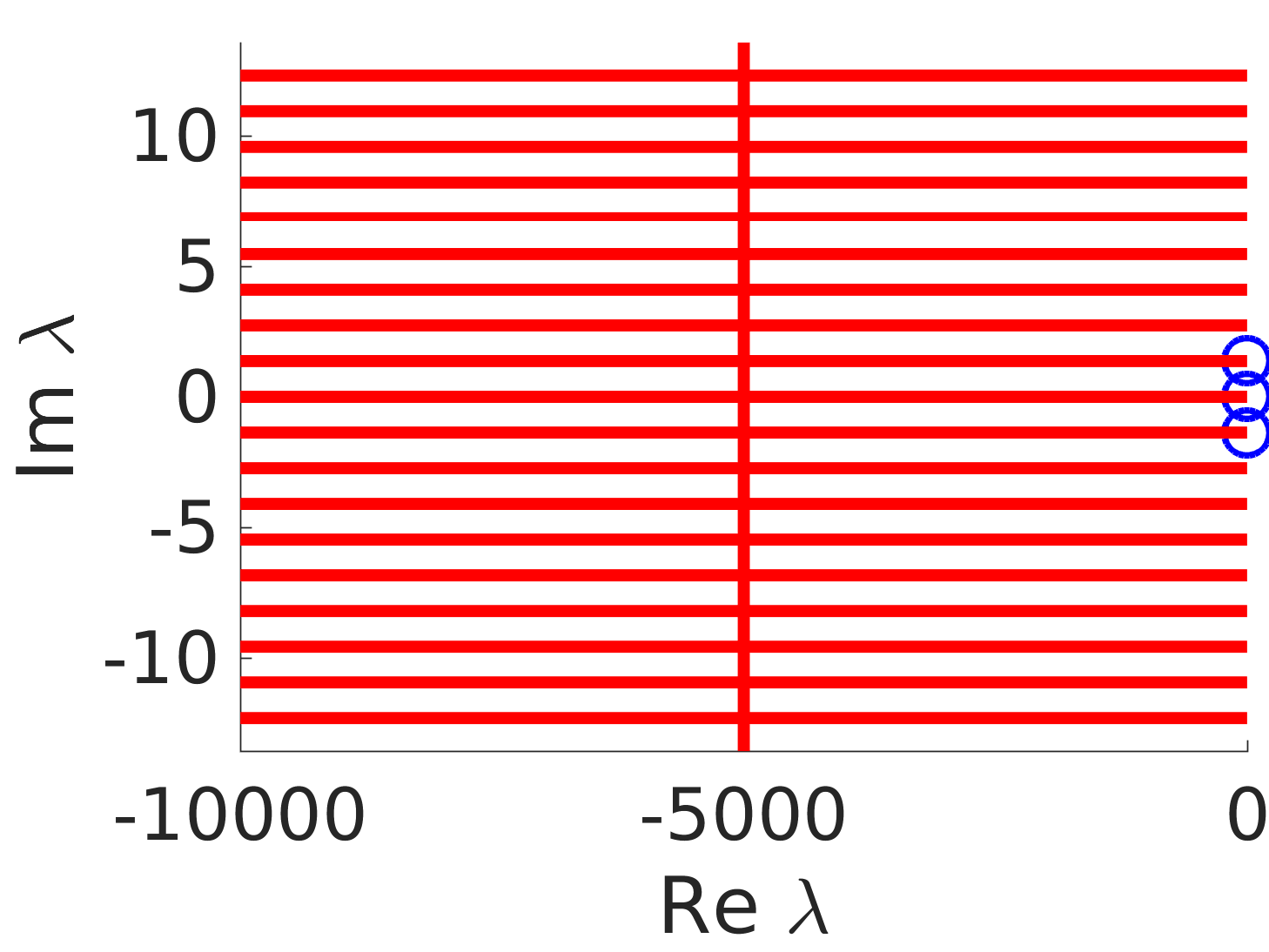}\label{fig:3.3a}}
  \subfigure[]{\includegraphics[height=4.2cm] {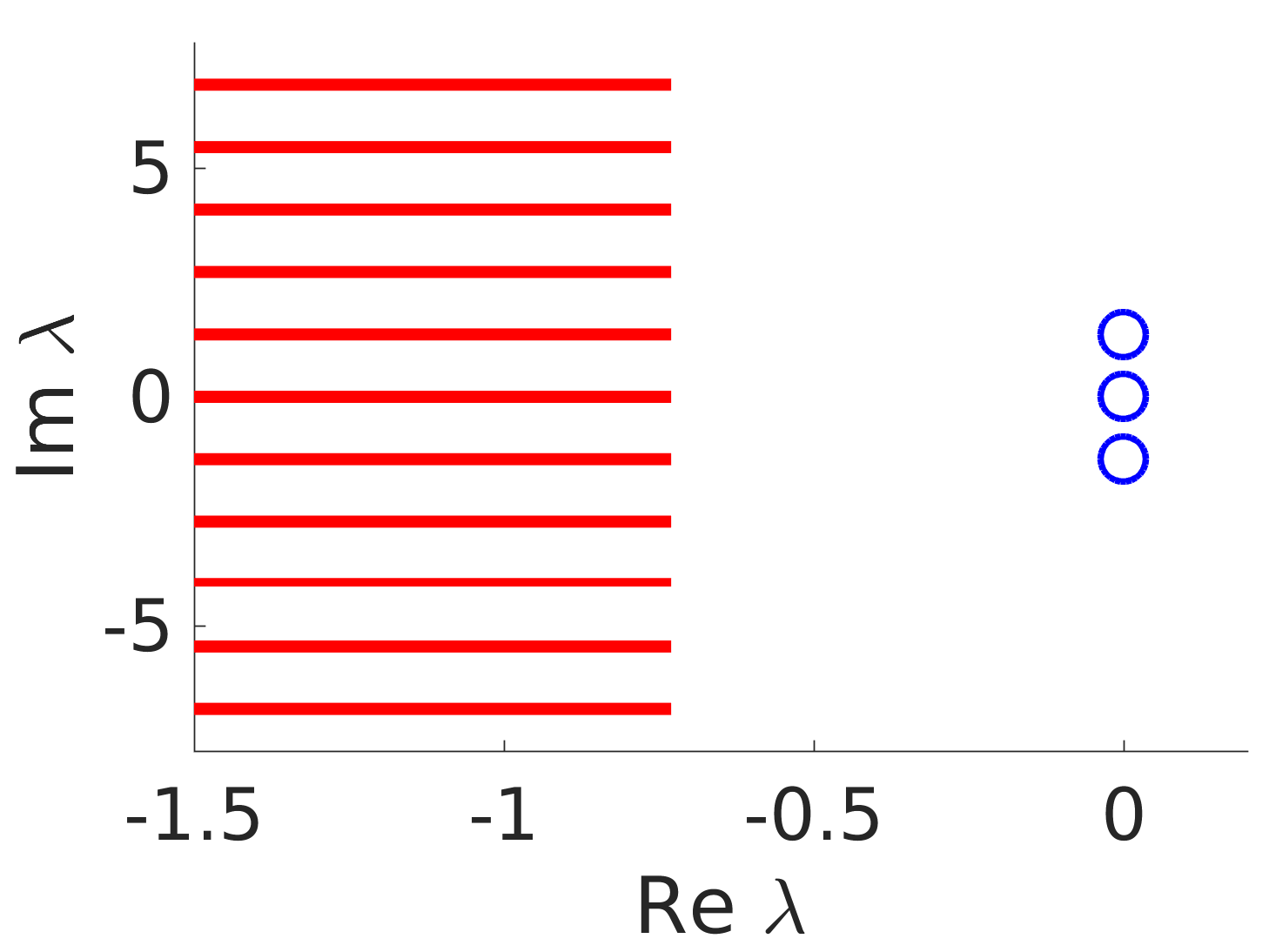}\label{fig:3.3b}}\\
  \subfigure[]{\includegraphics[height=4.2cm] {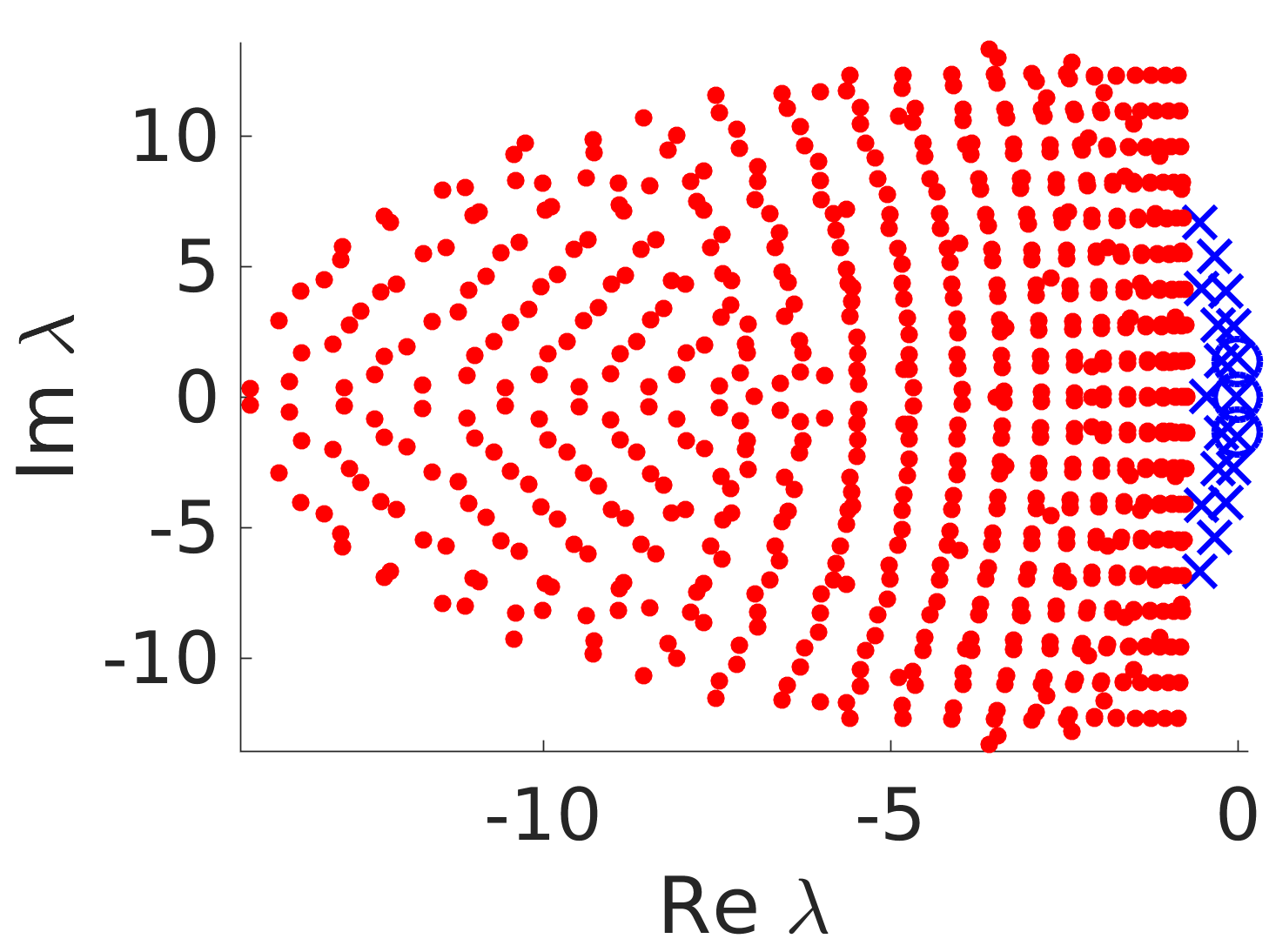} \label{fig:3.3c}}
  \subfigure[]{\includegraphics[height=4.2cm] {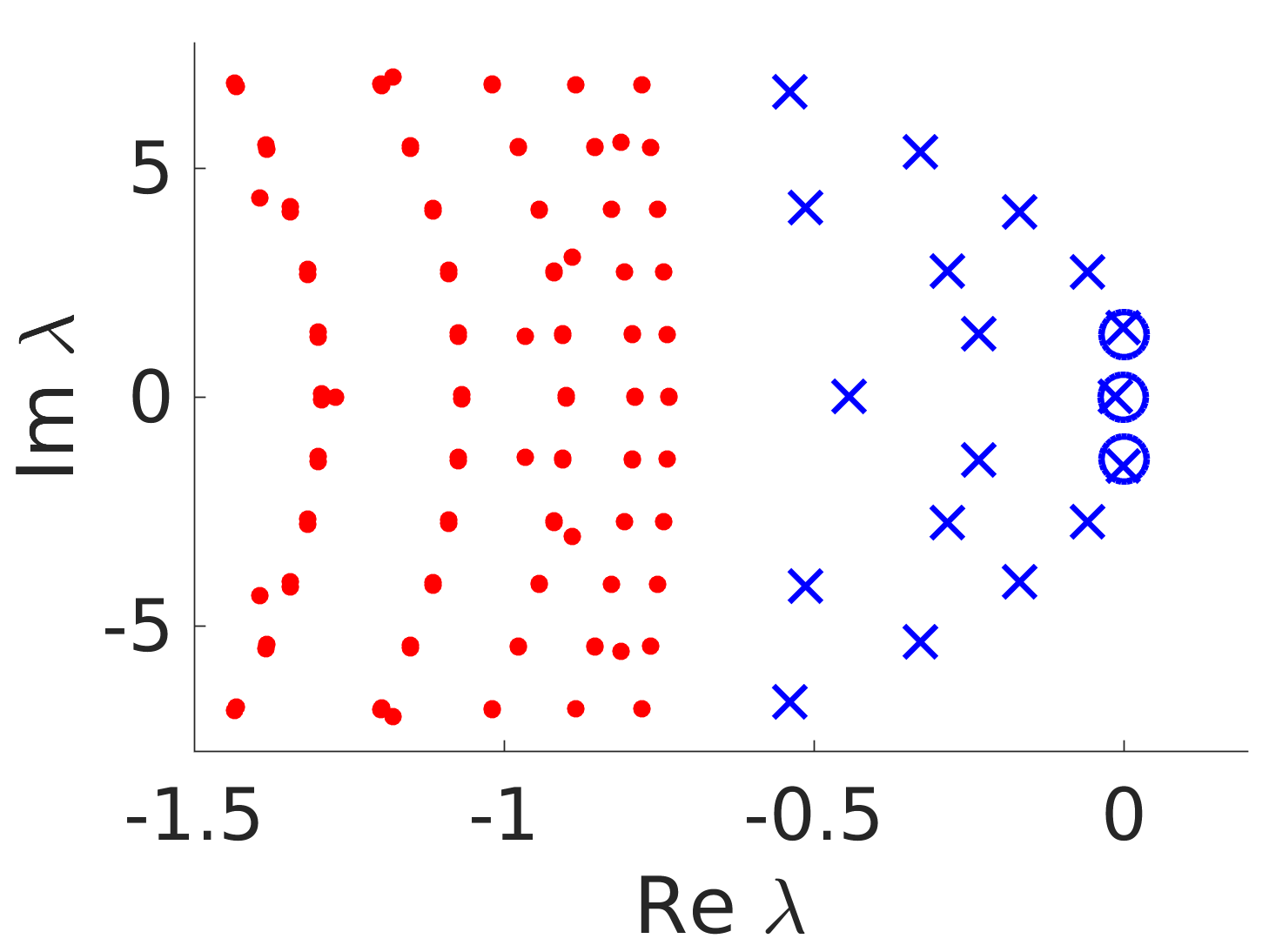}\label{fig:3.3d}}
  \caption{Subsets $\sigma_{\mathrm{disp}}(\mathcal{P})$ and $\sigma_{\mathrm{sym}}(\mathcal{P})$ of the spectrum for the cubic-quintic Ginzburg-Landau wave equation for $d=2$ with parameters \eqref{equ:3.41} (a),(b) and two different views of the numerical spectrum on a ball $B_R(0)$ with radius $R=20$ (c),(d).}
  \label{fig:3.3}
  \end{figure}

  Figure \ref{fig:3.3}(a) and (b) shows two different views for the part of the spectrum of the spinning solitons which is guaranteed by Proposition 
  \ref{prop:3.1} and \ref{prop:3.2}. It is subdivided into the symmetry set $\sigma_{\mathrm{sym}}(\PL)$ (blue circle), which is determined by Proposition \ref{prop:3.2} 
  and belongs to the point spectrum $\sigma_{\mathrm{pt}}(\PL)$, and the dispersion set $\sigma_{\mathrm{disp}}(\PL)$ (red lines), which is determined by Proposition \ref{prop:3.1} 
  and belongs to the essential spectrum $\sigma_{\mathrm{ess}}(\PL)$. In general, there may be further essential spectrum in $\sigma_{\mathrm{ess}}(\PL)\setminus\sigma_{\mathrm{disp}}(\PL)$ 
  and further isolated eigenvalues in $\sigma_{\mathrm{pt}}(\PL)\setminus \sigma_{\mathrm{sym}}(\PL) $. In fact, for the spinning solitons of the 
  cubic-quintic Ginzburg-Landau wave equation we find $18$ extra eigenvalues with negative real parts ($8$ complex conjugate pairs and $2$ purely real eigenvalues), 
  cf. Figure \ref{fig:3.3}(c),(d). These Figures show two different views for the numerical spectrum of the cubic-quintic Ginzburg-Landau wave equation on 
  the ball $B_R(0)$ with radius $R=20$ equipped with homogeneous Neumann boundary conditions. They consist of the approximations of the point spectrum 
  subdivided into the symmetry set (blue circle) and additional isolated eigenvalues (blue cross sign), and of the essential spectrum (red dots). 
  Three of these isolated eigenvalues are very close to the imaginary axis, see Figure \ref{fig:3.4}(c). Therefore, the spinning solitons seem to be 
  only weakly stable. Finally, the approximated eigenfunctions belonging to the eigenvalues $\lambda\approx 0$ and $\lambda\approx +i\sigma$ are shown 
  in Figure \ref{fig:3.4}(a) and (b). In particular, Figure \ref{fig:3.4}(a) is an approximation of the rotational term $\langle Sx,\nabla v_{\star}(x)\rangle$.

  \begin{figure}[H]
  \centering
  \subfigure[]{\includegraphics[height=4.0cm] {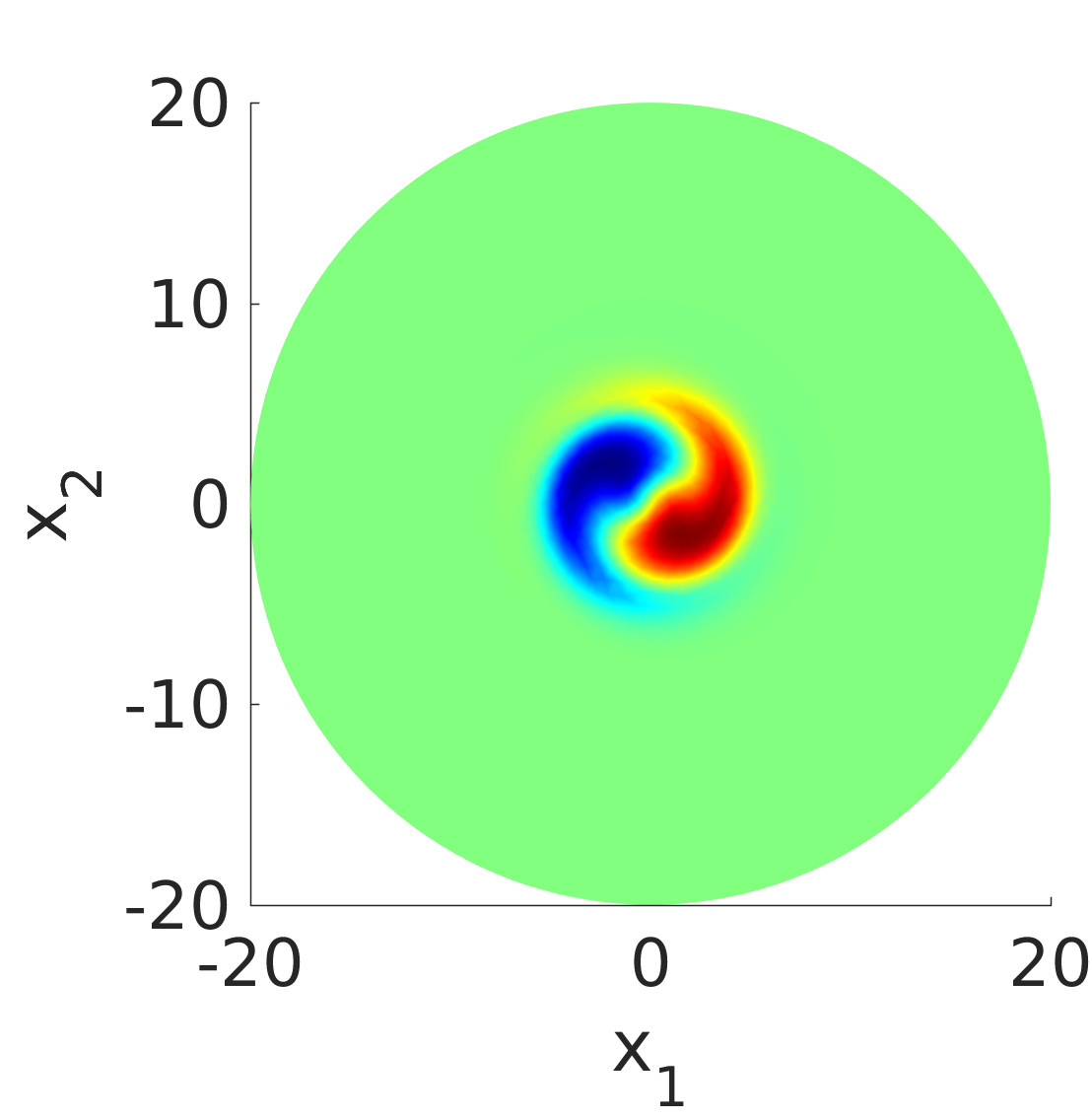}\label{fig:3.4a}}
  \subfigure[]{\includegraphics[height=4.0cm] {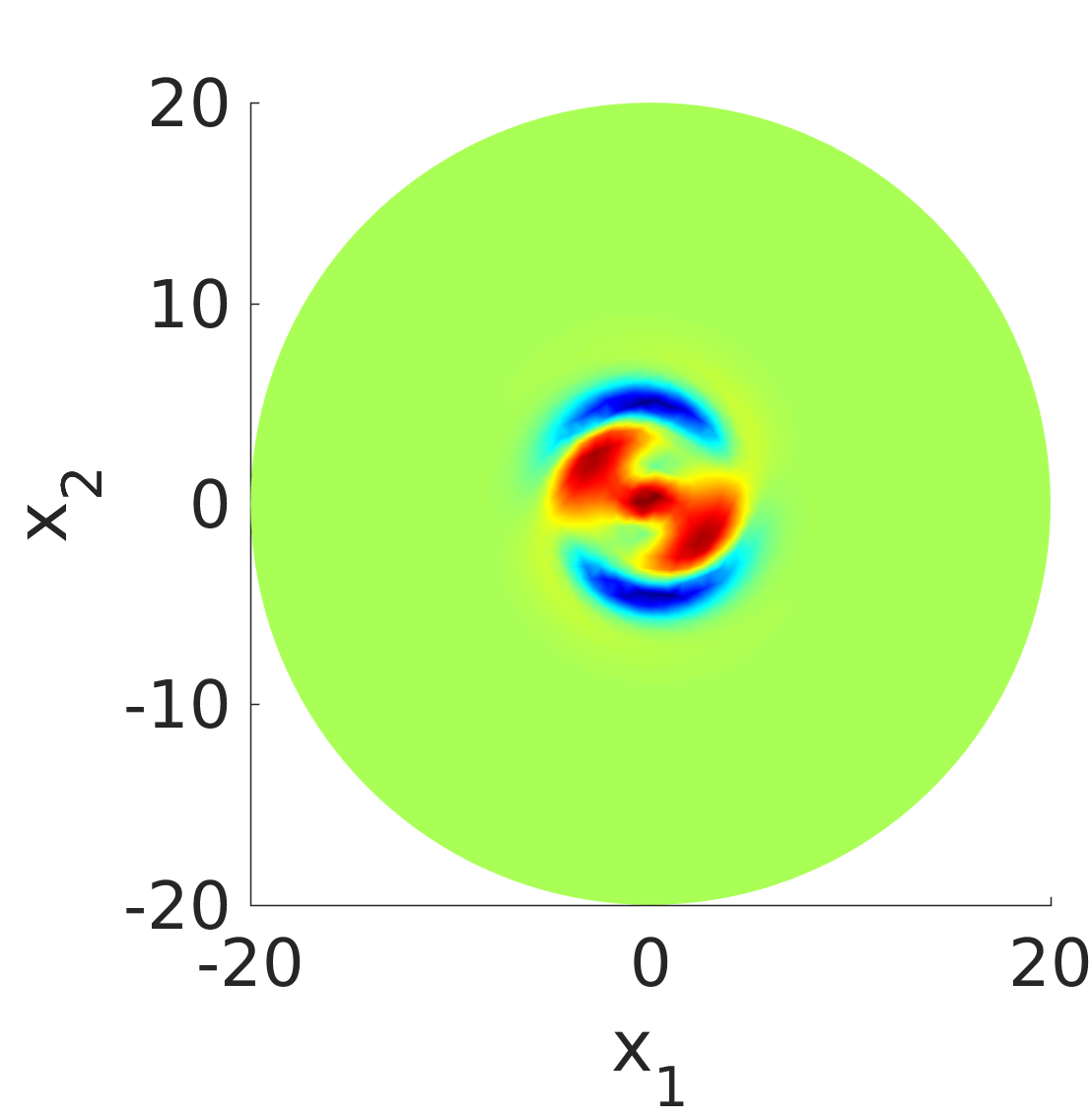}\label{fig:3.4b}}
  \subfigure[]{\includegraphics[height=4.0cm] {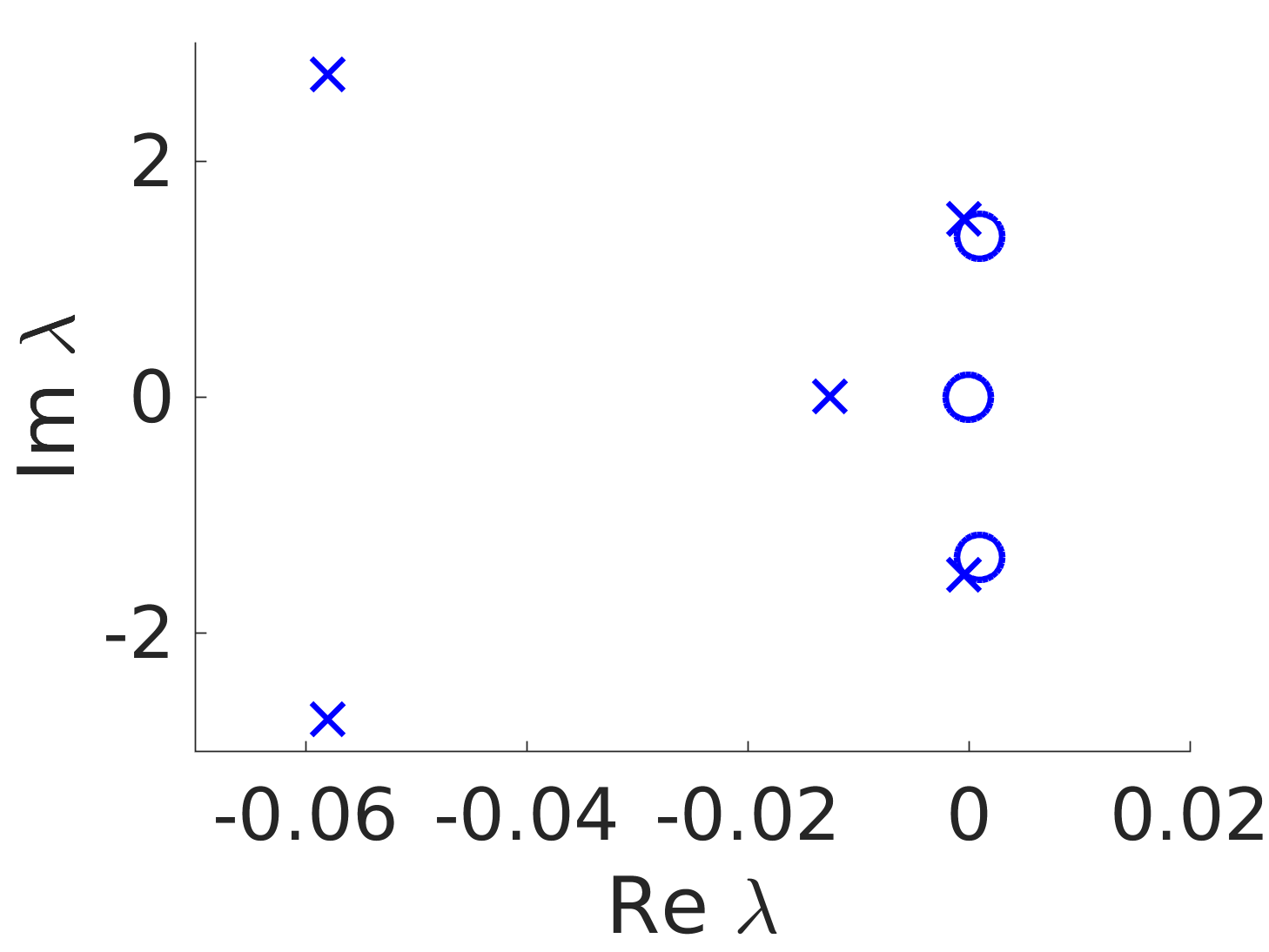}\label{fig:3.4c}}
  \caption{Eigenfunctions of the cubic-quintic Ginzburg-Landau wave equation for parameters \eqref{equ:3.41} belonging to the isolated eigenvalues $\lambda_1\approx 0$ (a) and $\lambda_2\approx i\sigma$ (b) and a zoom into the spectrum from Fig.\ref{fig:3.3}(c) in (c).}
  \label{fig:3.4}
  \end{figure}

\end{example}

\textbf{Acknowledgement.}
We gratefully acknowledge financial support by the Deutsche
      Forschungsgemeinschaft (DFG) through CRC 701 and CRC 1173.


\def\cprime{$'$}

\end{document}